\documentclass[a4paper,12pt,twoside]{article}
\usepackage{graphicx}\usepackage{a4wide}
\usepackage{amssymb,amsmath,amsthm}
\usepackage{yhmath} 
\usepackage{mathdots}
\usepackage[only,llbracket,rrbracket,interleave]{stmaryrd}

\usepackage{array}
\usepackage{multirow}
\usepackage{tabularx}
\usepackage{extarrows}
\usepackage{booktabs}

\usepackage{color,colordvi,graphicx,xcolor}
\usepackage{caption}
\usepackage{subcaption}

\usepackage[sort&compress]{natbib}

\usepackage{enumitem}

\newcommand{\bm}[1]{\text{\boldmath $#1$\unboldmath}}
\newcommand\bX{\bm{X}}
\newcommand\bU{\bm{U}}
\newcommand\bI{\bm{I}}
\newcommand\bK{\bm{K}}
\newcommand\bV{\bm{V}}
\newcommand\bSigma{\bm{\Sigma}}
\newcommand\bff{\bm{f}}

\newcommand\bz{\bm{z}}

\newcommand\bu{\bm{u}}

\newcommand\bn{\bm{n}}

\newcommand\bmu{\bm{\mu}}

\newcommand{\np}{{n_p}}

\newcommand\RR{\leavevmode\hbox{$\rm I\!R$}}
\newcommand{\TT}{\textsf{T}}
\newcommand\grad{\bm{\nabla}}

\newcommand{\lift}{\bm{r}_g}

\newcommand{\dbu}{\delta{\bu}}
\newcommand\bv{\bm{v}}
\newcommand\bw{\bm{w}}
\newcommand\dV{\,\textnormal{d}\Omega}

\newcommand\nodalu{{\textbf{u}}}
\newcommand\bNu{\bm{N}^{u}}
\newcommand\nodalp{{\textbf{p}}}
\newcommand\bNp{{N}^{p}}

\newcommand\bXu{\bm{X}_{\!u}}
\newcommand\bXp{\bm{X}_{\!p}}
\newcommand\bUu{\bm{U}_{\!\!u}}
\newcommand\bUp{\bm{U}_{\!\!p}}

\newcommand{\D}{\mathrm{D}}
\newcommand{\R}{\mathrm{R}}
\newcommand{\U}{\mathrm{U}}
\newcommand{\Uin}{\U_{\mathrm{in}}}
\newcommand{\Ujet}{\U_{\mathrm{jet}}}
\newcommand{\Ulid}{\U_{\mathrm{lid}}}
\newcommand{\UvecJet}{\mathbf{U}_{\!\mathrm{jet}}}
\newcommand{\UvecIn}{\mathbf{U}_{\!\mathrm{in}}}

\newcommand\nodaldu{{\bm{\delta} \!\textbf{u}}}
\newcommand\nodaldp{{\bm{\delta} \!\textbf{p}}}
\newcommand\nodalx{{\textbf{x}}}

\newcommand\nodalfu{{\textbf{f}}}

\newcommand\nodalfustar{{\textbf{f}_{u}^{\,\star}}}
\newcommand\nodalfpstar{{\textbf{f}_{p}^{\,\star}}}

\newcommand{\nS}{{\texttt{n}_{\texttt{S}}}}
\newcommand\ndofu{\texttt{n}^{\!u}}
\newcommand\ndofp{\texttt{n}^{p}}
\newcommand\ndofx{\texttt{n}^{x}}

\newcommand{\nk}{{\texttt{n}_{\texttt{POD}}}}
\newcommand{\nku}{{\texttt{n}_{\texttt{POD}}^{u}}}
\newcommand{\nkp}{{\texttt{n}_{\texttt{POD}}^{p}}}

\renewcommand\Re{{\text{Re}}}

\usepackage{scalerel,stackengine}
\stackMath
\newcommand\reallyhat[1]{%
\savestack{\tmpbox}{\stretchto{%
  \scaleto{%
    \scalerel*[\widthof{\ensuremath{#1}}]{\kern-.6pt\bigwedge\kern-.6pt}%
    {\rule[-\textheight/2]{1ex}{\textheight}}
  }{\textheight}%
}{0.5ex}}%
\stackon[1pt]{#1}{\tmpbox}%
}

\newenvironment{keywords}{\begin{quote}\emph{\textbf{Keywords:}}}{\end{quote}}

\theoremstyle{definition}

\newtheorem{remark}{Remark}

\usepackage{tikz}
\usetikzlibrary{fadings}
\usetikzlibrary{patterns}
\usetikzlibrary{shadows.blur}
\usetikzlibrary{shapes}

\graphicspath{ {figs/} }

\begin{document}
\title{Data augmentation for the POD formulation of the parametric laminar incompressible Navier-Stokes equations}

\author{
\renewcommand{\thefootnote}{\arabic{footnote}}
			Alba Muix\'i\footnotemark[1]\textsuperscript{ \ ,}* , \
			Sergio Zlotnik\footnotemark[1]\textsuperscript{ \ ,}\footnotemark[2] , \\ 
\renewcommand{\thefootnote}{\arabic{footnote}}
			Matteo Giacomini\footnotemark[1]\textsuperscript{ \ ,}\footnotemark[2]  \  and
			Pedro D\'iez\footnotemark[1]\textsuperscript{ \ ,}\footnotemark[2]
}

\date{}
\maketitle

\renewcommand{\thefootnote}{\arabic{footnote}}

\footnotetext[1]{Laboratori de C\`alcul Num\`eric (LaC\`aN), ETS de Ingenier\'ia de Caminos, Canales y Puertos, Universitat Polit\`ecnica de Catalunya, Barcelona, Spain.}
\footnotetext[2]{Centre Internacional de M\`etodes Num\`erics en Enginyeria (CIMNE), Barcelona, Spain.
\vspace{5pt}\\
* Corresponding author: Alba Muix\'i \textit{E-mail:} \texttt{alba.muixi@upc.edu}
}

\begin{abstract}
A posteriori reduced-order models (ROM), e.g. based on proper orthogonal decomposition (POD), are essential to affordably tackle realistic parametric problems. They rely on a  trustful training set, that is a family of full-order solutions (snapshots) representative of all possible outcomes of the parametric problem. Having such a rich collection of snapshots is not, in many cases, computationally viable. A strategy for data augmentation, designed for parametric laminar incompressible flows, is proposed to enrich poorly populated training sets. The goal is to include in the new, artificial snapshots emerging features, not present in the original basis, that do enhance the quality of the reduced basis (RB) constructed using POD dimensionality reduction. The methodologies devised are based on exploiting basic physical principles, such as mass and momentum conservation,  to construct physically-relevant, artificial snapshots at a fraction of the cost of additional full-order solutions. Interestingly,  the numerical results show that the ideas exploiting only mass conservation (i.e., incompressibility) are not producing significant added value with respect to the standard linear combinations of snapshots.  Conversely, accounting for the linearized momentum balance via the Oseen equation does improve the quality of the resulting approximation and therefore is an effective data augmentation strategy in the framework of viscous incompressible laminar flows. Numerical experiments of parametric flow problems, in two and three dimensions,  at low and moderate values of the Reynolds number are presented to showcase the superior performance of the data-enriched POD-RB with respect to the standard ROM in terms of both accuracy and efficiency.
\end{abstract}


\section{Introduction}

Reduced-Order Models (ROM) are commonly employed to construct affordable solutions of multi-query, parametric problems \citep{Chinesta:2017:ECM,Willcox-PWG-18}.
Although well established, the techniques to devise ROMs might face computational bottlenecks, when the underlying physical models are nonlinear \citep{Farhat-CBF-11,Farhat-CFCA-13,nguyen2023efficient},  the space of parameters is high dimensional \citep{Willcox-BWG-08,Constantine-15} and the cost of each simulation is expensive, like in computational fluid dynamics applications \cite{Sherwin-MMTMJTPS-21,Lohner-LOMFD-21}.
Indeed, the issues mentioned above are particularly critical when parametric incompressible flows are considered, from laminar \citep{Stabile-SR-18,TsiolakisPGD2020}, to turbulent \citep{Stabile-HSMR-20,Iliescu-APSRIN-21,Tsiolakis-TGSOH-22} regimes,  in the context of parameterized geometries \citep{Rozza-BFIMQRS-16,Giacomini-GBSH-21} or uncertain inputs \citep{Nobile-MN-18,Tamellini-PTPBSD-23}.
While such challenges are particularly relevant in transient and turbulent flows, many conceptual difficulties already arise in steady and laminar flows. Specifically, this work focuses on the parametric incompressible Navier-Stokes equations at low and moderate Reynolds numbers.

To construct a surrogate model for such parametric flow problems, a posteriori ROMs rely on collecting a family of snapshots (corresponding to different instances of the parameters) and using them as a basis to describe a reduced space where the solution is sought \citep{burkardt2006podcvt,Patera-Rozza:07,rozza2008reduced,nguyen2008efficient,wang2012podturbulent}. 
It is worth noticing that the set of snapshots, the training set, has to be representative of the full space of solutions, and this generally entails a large number of high-fidelity solutions, and a high computational cost.  Possible redundancies in the training set typically yield ill-conditioned reduced-order problems. These redundancies are however readily suppressed using Proper Orthogonal Decomposition (POD). 

Of course, the outcome of these operations depends on the quality of the data in the original training set.
It is well known that the selection of \emph{optimal} and \emph{adaptive} sampling procedures is crucial to minimize the number of required full-order computations during the construction of a ROM \cite{Veroy-VP-05,Willcox-BWG-08,Hesthaven-HSZ-14,Breitkopf-PBBVZ-20}.
Nonetheless, in high-dimensional parametric spaces, this operation is not trivial, see, for instance, the review paper \cite{Vono-VDC-22}.
More recently, alternative approaches to make the construction of ROMs computationally affordable relied on multi-fidelity strategies mixing simulations on different levels of mesh refinement in the physical and parametric space \citep{Willcox-NW-14,Breitkopf-XZBVZ-18,Jakeman-JFETGA-22} and domain decomposition techniques to couple local ROMs computed on smaller subdomains \citep{Barnett:2022:Sandia,Hesthaven-DH-23,MG-DEG-24}.

Although effective, all the above-mentioned strategies entail the solution of a certain number of full-order problems to populate the training set. 
In particular, when little to no prior knowledge of the solution is available, two main families of approaches are suitable to treat such problems. On the one hand,  the paradigm of a priori ROM \cite{Chinesta-review-13,MG-GBSH-21} circumvents the sampling issue by constructing a reduced model without preliminarily computing any snapshot. On the other hand,  data augmentation techniques \cite{diez2021nonlinear,Nguyen-Preprint-24} aim to emancipate the quality of the a posteriori ROM approximation from the richness of the initial dataset, incorporating new information in uniformly-sampled, and possibly poorly populated, existing training sets.

This work proposes to reduce the computational investment needed to generate the training set by producing artificial snapshots and augmenting data in the context of parametric viscous laminar incompressible flows.
The idea stems from the work in \citep{diez2021nonlinear}, where a strategy to expand the training set by adding new elements, artificially generated from the original snapshots with simple algebraic operations,  was proposed. Of course,  these artificial snapshots are meaningful if they enrich the approximation space, providing better approximations to new solutions of the parametric problem. 
It is worth recalling that augmenting the training set is never counterproductive (POD suppresses redundancies) and it can be extremely useful if the new elements bring pieces of information that are absent in the original set and pertinent to capture features of new parametric solutions.
The present contribution bypasses the simple framework of linear convection-diffusion equation to treat parametric viscous laminar incompressible flows, modeled by the steady incompressible Navier-Stokes equations. This family of problems entails two new challenges not accounted for in previous data augmentation approaches: (i) the incompressibility constraint and (ii) the nonlinearity of the convection term.
To create new training solutions combining the original snapshots,  this work proposes a set of new data augmentation strategies that exploit the knowledge of the physics underlying the problem under analysis. 
The resulting physics-informed data augmentation enforces the mass and momentum conservation equations to construct artificial snapshots, at a fraction of the cost of computing additional full-order solutions. These yield new, relevant information not present in the original dataset, while guaranteeing that the balance laws for the flow system are fulfilled, for any combination of snapshots considered.

The remainder of the paper is structured as follows. Section \ref{sec:NS} describes the problem statement, introduces the notation and the standard strategy to solve the full-order model. Then, in Section \ref{sec:POD}, the basics of POD are briefly recollected, together with its application to the reduced basis (RB) formulation of mixed nonlinear problems. Section \ref{sec:data_augmentation} devises the proposed methodologies to augment the data and enrich the training set, accounting for the physical knowledge of the problem at hand.
Numerical examples, presented in Section \ref{sec:NumEx}, demonstrate the suitability of the discussed approaches to augment the training set with physically-relevant information, enhancing the performance of the POD-RB approximation.
Finally, a discussion on the accuracy and efficiency of the proposed methodology is presented in Section \ref{sec:Discussion} and the concluding remarks are collected in Section \ref{sec:Conclusion}.

\section{Full-order model: steady incompressible Navier-Stokes equations}
\label{sec:NS}
In this Section, the full-order model of the problem under analysis, the steady laminar incompressible Navier-Stokes equations, is recalled and the high-fidelity solver employed for the computation of the snapshots is introduced.

\subsection{Problem statement}
The problem under consideration consists in finding a velocity field $\bu$ and a pressure field $p$ taking values in $\Omega\subset\RR^d$ (the space dimension $d$ is equal to 2 or 3),  and such that 
\begin{equation}\label{eq:NavierStokes}
\left\{
\begin{aligned}
-\grad \cdot \left( \nu \grad \bu - p \bI \right) +  \left( \bu \cdot \grad \right) \bu &= \bm{0}  &\text{in } \Omega, \\
\grad \cdot \bu &= 0  & \text{in } \Omega, \\
\bu &= \bu_D & \text{on } \Gamma_D, \\
\left(  \nu \grad \bu - p \bI \right) \bn &= \bm{t} & \text{on } \Gamma_N,
\end{aligned}
\right.
\end{equation}
where $\nu>0$ is the kinematic viscosity,  $\bu_D$ the velocity prescribed on $\Gamma_D$, $\bm{t}$ the pseudo-traction prescribed on $\Gamma_N$, 
$\bI$ stands for the $d\times d$ identity matrix, and $\bn$ is the outward unit normal to the boundary. 
The two parts $\Gamma_D$ and $\Gamma_N$ of the boundary $\partial\Omega$, where Dirichlet and Neumann boundary conditions are enforced, are such that
$\partial\Omega = {\bar\Gamma_D} \cup {\bar\Gamma_N}$, and $\Gamma_D\cap\Gamma_N = \emptyset$.

Problem \eqref{eq:NavierStokes} is considered to be parametric in the sense that the input data (and therefore the solution $(\bu , p)$) depend on $\np$ parameters collected in a vector $\bmu\in\RR^\np$. These parameters may affect material properties of the fluid ($\nu$), working conditions ($\bu_D$, $\bm{t}$), the geometry of the domain or the parts of the boundary where the boundary conditions are enforced.

The corresponding weak form is derived using weighted residual technique and reads: 
find $(\bu , p) \in [H^1(\Omega)]^d \times L^2(\Omega)$ such that $\bu = \bu_D$ on $\Gamma_D$ and
\begin{equation}\label{eq:WeakFormNS}
\left\{
\begin{aligned}
d(\bu,\bv) + c(\bu,\bu,\bv) - b(\bv,  p) &= \int_{\Gamma_N} \bm{t} \cdot \bv \,d\Gamma, \\
b(\bu,q) & = 0 ,
\end{aligned}
\right.
\end{equation}
for every $\bv \in [H^1_{0}(\Omega)]^d  = \{\bv \in [H^1(\Omega)]^d : \bv|_{\Gamma_D} = \bm{0} \}$ and for every $q \in L^2(\Omega)$. The forms in \eqref{eq:WeakFormNS} are
\begin{equation}\label{eq:Forms}
\begin{aligned}
d(\bu,\bv) &= \int_\Omega \nu \, \grad \bu : \grad \bv  \dV, \\
c(\bw,\bu,\bv) &= \int_\Omega \left( \bw \cdot  \grad \bu \right) \cdot \bv \dV, \\ 
b(\bv,  q) &= \int_\Omega q\, \grad \cdot \bv \dV.
\end{aligned}
\end{equation}

In practice, one possibility to deal with inhomogeneous Dirichlet conditions ($\bu_D \neq \bm{0}$) is taking $\bu=\bu_{0}+\lift$,  $\lift$ being a lift function such that $\lift |_{\Gamma_D} = \bu_D$. Thus, the actual unknown is such that $\bu_{0}|_{\Gamma_D} = \bm{0}$ and belongs to the same space as $\bv$. 

Note that problem \eqref{eq:NavierStokes} is nonlinear and therefore is to be solved iteratively. Next Section briefly describes the formulation of the Newton-Raphson method employed within the high-fidelity solver in this work.

\subsection{Newton-Raphson iteration}
An iterative scheme is devised to numerically solve problem  \eqref{eq:NavierStokes} (or the equivalent weak form \eqref{eq:WeakFormNS}). Given some approximation $(\bu^k, p^k)$, for $k = 0,1,\dots$ ($k$ is the iteration counter), the subsequent iteration $(\bu^{k+1}, p^{k+1})$ is characterized by increments $(\dbu,\delta p)$ such that $\bu^{k+1}= \bu^{k} + \dbu$ and $p^{k+1}=p^{k}+\delta p$. 

Thus, introducing the $\bu= \bu^{k} + \dbu$ and $p=p^{k}+\delta p$ in  \eqref{eq:WeakFormNS} results in the following linearized equation for the unknown $(\bu,p)$,
\begin{equation}\label{eq:WeakFormIter}
\left\{
\begin{aligned}
d(\dbu,\bv) + c(\bu^{k},\dbu,\bv) + c(\dbu,\bu^{k},\bv) - b(\bv,  \delta p) &= \int_{\Gamma_N} \bm{t} \cdot \bv \,d\Gamma  \\
& \,\,\,\,\, -d(\bu^{k},\bv) \\
& \,\,\,\,\, - c(\bu^{k},\bu^{k},\bv) \\
& \,\,\,\,\, + b(\bv, p^{k}),\\
b(\dbu,q) & = -b(\bu^{k},q).
\end{aligned}
\right.
\end{equation}

Note that the term $c(\dbu,\dbu,\bv)$ features a second-order increment, thus it is considered to be negligible and drops out in \eqref{eq:WeakFormIter}.

It is assumed that $\bu^{k}$ fulfils the Dirichlet boundary conditions, therefore $\dbu$ is enforced to satisfy their homogeneous counterpart.

The linear problem \eqref{eq:WeakFormIter} characterising every iteration is discretized with a standard finite element approach. 
Namely, the velocity field $\bu^{k}$ is represented by the vector of nodal values $\nodalu^{k}\in \RR^{\ndofu}$, and the pressure field $p^{k}$ is represented, in its discrete form, by the vector $\nodalp^{k}\in \RR^{\ndofp}$ ($\ndofu$ and $\ndofp$ are the number of degrees of freedom in the velocity and pressure discretizations), such that
\begin{equation}
\bu^{k}=\sum_{i=1}^{\ndofu} [\nodalu^{k}]_{i} \bNu_{i} 
\quad \text{ and } \quad
p^{k}=\sum_{i=1}^{\ndofp} [\nodalp^{k}]_{i} \bNp_{i} ,
\end{equation}
$\bNu_{i}$ and $\bNp_{i}$, being the shape functions corresponding to the $i$-th degree of freedom of velocity and pressure, respectively, while $[\nodalu^{k}]_{i} $ and $[\nodalp^{k}]_{i}$ denote the $i$-th component of the vectors of nodal values $\nodalu^{k}$ and $\nodalp^{k}$.
Henceforth,  a finite element pair for velocity and pressure fulfilling the inf-sup or LBB condition, e.g., the Taylor-Hood elements, is assumed \citep{donea2003finite,quarteronivalli2008}.
Analogously, the vectors of nodal values describing $\dbu$ and $\delta p$ read $\nodaldu$ and $\nodaldp$. We are aware the difference of notation between the fields and the vectors of nodal values representing them is somehow subtle. Nevertheless, we firmly believe the reader will distinguish them easily due to the context.

The matrices representing the discrete counterparts of the operators introduced in \eqref{eq:Forms}, in the discrete spaces spanned by the finite element basis functions $\bNu_{i}$ (for $i=1,2,\dots,\ndofu$) and $\bNp_{i}$ (for $i=1,2,\dots,\ndofp$), are
\begin{equation*}
\begin{aligned}
\bm{D}\in\RR^{\ndofu\times\ndofu} & \text{ with generic term } &[\bm{D}]_{ij} &=d(\bNu_{i},\bNu_{j}),\\
\bm{C}_1(\bu^k)\in\RR^{\ndofu\times\ndofu} & \text{ with generic term } & [\bm{C}_1(\bu^k)]_{ij}&=c(\bu^k, \bNu_{j},\bNu_{i}),\\
\bm{C}_2(\bu^k)\in\RR^{\ndofu\times\ndofu} & \text{ with generic term } & [\bm{C}_2(\bu^k)]_{ij}&=c(\bNu_{j},\bu^k,\bNu_{i}),\\
\bm{B}\in\RR^{\ndofp\times\ndofu} & \text{ with generic term } & [\bm{B}]_{ij}&=b(\bNu_{j},\bNp_{i}).
\end{aligned}
\end{equation*}
Thus, the discrete version of \eqref{eq:WeakFormIter} results in the following linear algebraic system:
\begin{equation}\label{eq:LinearSystemIter}
\begin{bmatrix}
\bm{D}+ \bm{C}_1(\bu^k) +\bm{C}_2(\bu^k)& \bm{B}^\TT \\
\bm{B} & \bm{0} 
\end{bmatrix}
\begin{bmatrix}
\nodaldu  \\
\nodaldp
\end{bmatrix} = 
\begin{bmatrix}
\nodalfu - \bm{D}\nodalu^k - \bm{C}_1(\bu^k) \nodalu^k  - \bm{B}^\TT \nodalp^k  \\
- \bm{B} \nodalu^k 
\end{bmatrix},
\end{equation} 
where $\nodalfu$ is the discrete version of the Neumann term in  \eqref{eq:WeakFormIter}, that is, 
$$
[\nodalfu \,]_{i}=  \int_{\Gamma_N} \bm{t} \cdot \bNu_i \,d\Gamma \, ,  \quad \text{ for } i=1,2,\dots,\ndofu \, .
$$

Thus, system \eqref{eq:LinearSystemIter} is to be solved at each iteration. Both the matrix and the right-hand-side term depend on the input parameters of the system (the parameters $\bmu$ mentioned above) and on the previous iteration $(\bu^k, p^k)$. The need of solving many systems corresponding to slightly different parametric entries suggests the use of ROMs.

\section{Reduced-order model: Proper Orthogonal Decomposition}
\label{sec:POD}

In this Section, the main ideas of the standard Proper Orthogonal Decomposition (POD) formulation are briefly summarized, with special emphasis towards the construction of a stable reduced-order formulation, POD-RB, for incompressible flows \citep{Veroy-VP-05,Veroy-RV-07,ballarin2015supremizer,Stabile-SR-18}.

\subsection{POD-RB for linear systems of equations}
\label{sec:POD-linearSys}

Consider the general case in which the parametric unknown $\nodalx(\bmu)$ is a vector of dimension $\ndofx$ obtained as the solution of the linear system of equations
\begin{equation}\label{eq:fullOrder}
\bK(\bmu) 	\nodalx(\bmu) = \bff(\bmu),
\end{equation}
where, with respect to the previous Section, $\nodalx$ plays the role of $\nodaldu$ and $\nodaldp$.

Let $\nS$ be the number of snapshots $\nodalx_i$, for $i = 1,\dots,\nS$ corresponding to high-fidelity solutions of the full-order problem \eqref{eq:fullOrder} for $\nS$ values of the input parameters $\bmu$.
The snapshots are centred and collected in the matrix $\bX$ such that 
\begin{equation}\label{eq:matrixX}
\bX = \left[ \widehat{\nodalx}_1\,\,  \widehat{\nodalx}_2\,\,  \dots \,\, \widehat{\nodalx}_{\nS} \right], \text{ where }
\widehat{\nodalx}_i ={\nodalx}_i - \bar{\nodalx} \text{ and } 
\bar{\nodalx} =  \dfrac{1}{\nS} \sum_{i=1}^{\nS} {\nodalx}_{i}.
\end{equation}

In some cases, scaling or normalizing the snapshots in $\bX$ is recommended to avoid overrepresentation of some of them. In the presented examples in Section \ref{sec:NumEx}, only pressure snapshots in Section \ref{sc:microswimmer} are scaled since they present significantly different orders of magnitude.

The POD, which is based on the Singular Value Decomposition (SVD), is used to eliminate redundancies in $\bX$. The SVD of $\bX$ reads
\begin{equation}
\bX = \bU \bSigma \bV^{\TT},
\end{equation}
where $\bSigma\in\RR^{\ndofx\times\nS}$ is a diagonal matrix containing the singular values of $\bX$ in descending order, $\sigma_1\ge\sigma_2\ge\dots\ge\sigma_\nS\ge0$.

The first $\nS$ columns of $\bU$ define an orthonormal basis for the linear subspace generated by the snapshots. The more relevant information is contained in the first ones, which correspond to larger singular values. In practice, selecting a tolerance $\varepsilon$, the first $\nk$ columns of $\bU$ are kept, $\nk$ being  such that
\begin{equation}\label{eq:criterion}
\sum_{i=1}^\nk \sigma_i \ge (1-\varepsilon) \sum_{i=1}^{\nS} \sigma_i.
\end{equation}

The POD consists in representing the unknown $\nodalx$ in terms of only these $\nk$ modes ($\nk$ is expected to be much lower than $\nS$, and $\nS$ is typically taken much lower than the dimension of the full-order system). Namely,
\begin{equation}\label{eq:POD_approx}
\nodalx \simeq \bar\nodalx  + \bU^\star \bz =  \bar\nodalx+ \sum_{i=1}^\nk [\bU]_i z_i,
\end{equation}
where $\bz=[z_1 \, z_2 \dots z_\nS]^\top$ is the vector of the POD unknowns, and $\bU^\star$ is the matrix containing the first $\nk$ columns of $\bU$. Thus, following the RB rationale, a Galerkin projection of the full-order system \eqref{eq:fullOrder} yields the POD-RB system
\begin{equation}\label{eq:reducedOrder}
[\bU^{\star \top}  \bK(\bmu) \bU^\star ] \, \bz = \bU^{\star \top} \bff(\bmu) - \bU^{\star \top}  \bK(\bmu) \, \bar\nodalx ,
\end{equation}
of reduced size $\nk \ll \nS$.

\subsection{POD-RB for Navier-Stokes iterations}
\label{sec:PODNS}
As mentioned in Section \ref{sec:POD-linearSys}, $\nS$ snapshots are collected as solutions of the discrete version of the full-order problem \eqref{eq:WeakFormNS} corresponding to different values of the input parameters $\bmu$, that is the solution of \eqref{eq:LinearSystemIter} upon convergence of the Newton-Raphson iterations.

Two snapshot matrices  $\bXu$ and $\bXp$ are thus collected separately for velocity and pressure. Following the POD strategy described above, the SVD applied to $\bXu$ and $\bXp$ provides truncated bases $\bUu^{\star}$ and $\bUp^{\star}$, being $\nku$ and $\nkp$ the number of terms kept for each approximation. Note that the tolerance $\varepsilon$ used in \eqref{eq:criterion} to set the number of terms may be different for velocity and pressure.

It follows that the unknowns $\nodaldu$ and $\nodaldp$ are approximated as 
\begin{equation}\label{eq:approxDuDp}
\nodaldu \simeq \bUu^{\star} \bm{z}_u \quad \text{ and } \quad \nodaldp \simeq \bUp^{\star} \bm{z}_p,
\end{equation}
where $\bm{z}_u$ and $\bm{z}_p$ are the reduced unknowns of size $\nku$ and $\nkp$, respectively.
In addition, the reduced version of the problem to be solved at each Newton-Raphson iteration is obtained pre-multiplying the equations in system \eqref{eq:LinearSystemIter} by $\bUu^{\star \, \top}$ and $\bUp^{\star \, \top}$, yielding
\begin{subequations}\label{eq:ReducedLinearSystemIter} 
	\begin{equation}\label{eq:ReducedLinearSystem} 
	\begin{bmatrix}
	\bK^\star & \bm{B}^{\star \top} \\
	\bm{B}^{\star} & \bm{0} 
	\end{bmatrix}
	\begin{bmatrix}
	\bm{z}_u \\
	\bm{z}_p
	\end{bmatrix} = 
	\begin{bmatrix}
	\nodalfustar  \\
	\nodalfpstar
	\end{bmatrix},
	\end{equation} 
where the following matrices and vectors are introduced
	\begin{align}\label{eq:MatricesVectors}
	\bK^\star &= \bUu^{\star \, \top} \left[ \bm{D}+ \bm{C}_1(\bu^k) +\bm{C}_2(\bu^k)  \right] \bUu^{\star} \,,
	\\
	\bm{B}^\star &=  \bUp^{\star \, \top} 	\bm{B} \bUu^{\star} \, ,
	\\
	\nodalfustar &= \bUu^{\star \, \top} \left[  \nodalfu - \bm{D}\nodalu^k - \bm{C}_1(\bu^k) \nodalu^k  - \bm{B}^\TT \nodalp^k \right] \, ,
	\\
	\nodalfpstar &=  -\bUp^{\star \, \top} 	\bm{B} \nodalu^k \, .
	\end{align}
\end{subequations}

Note that the reduced system \eqref{eq:ReducedLinearSystemIter} maintains the saddle-point structure of the full-order approximation \eqref{eq:LinearSystemIter}.
It is well-known that, although the snapshots are computed using an LBB-compliant discretization, it is not guaranteed that the stability properties of the full-order system are preserved at the reduced level \citep{Veroy-VP-05,Veroy-RV-07,Veroy-GV-12}.

A possible workaround for this issue relies on avoiding to reconstruct the pressure in the reduced space. Hence,  the POD-RB formulation of the saddle-point problem \eqref{eq:ReducedLinearSystemIter} is only employed to construct a reduced approximation of velocity (see, e.g., \cite{ito1998rom,Veroy-VP-05,grassle2019pod}).  Nonetheless, this approach poses a series of difficulties in problems requiring the evaluation of aerodynamic forces \citep{Noack-NPM-05}.

Moreover, it is worth recalling that, due to the linearity of the mass conservation equation,  if all the elements in the velocity basis are weakly incompressible,  any linear combination of such terms automatically verifies the divergence-free condition on velocity. Since the solenoidal nature of the velocity is preserved, it is no longer required to introduce a Lagrange multiplier (i.e., pressure) to enforce the incompressibility constraint \citep{Adeline-MFPH-10}.
Whilst it is thus straightforward to enforce incompressibility at the weak level when homogeneous Dirichlet conditions are imposed on all the boundaries, the cases featuring inhomogeneous Dirichlet conditions (possibly depending upon the parameters $\bmu$) or Neumann conditions are not so trivial.  Under these assumptions, a linear combination of the snapshots, complemented with the  Dirichlet boundary conditions, is not automatically solenoidal. 
Therefore, in order to enforce that the reduced-order approximation of the velocity is divergence-free, a Lagrange multiplier is indeed required also at the reduced level. This entails the need to maintain the pressure, interpreted as the Lagrange multiplier to impose the incompressibility constraint \citep{donea2003finite},  in equation \eqref{eq:ReducedLinearSystemIter}, yielding a mixed formulation in the reduced space.
Thus, despite its attractiveness, the use of solenoidal velocity snapshots to avoid solving a saddle-point problem for the velocity-pressure pair in the reduced space is not a viable option if the boundary conditions are not homogeneous.

\subsection{Stability and pressure reconstruction}
\label{sec:Stability}

It is well-known that the pressure term is needed in many applications \citep{Noack-NPM-05} and many quantities of engineering interest (e.g., aerodynamic forces, pressure drops, $\ldots$) require an assessment of the pressure field.

Nonetheless, as previously mentioned, the inf-sup stability of the full-order discretization does not imply the inf-sup stability of the reduced system.
To achieve stable reduced approximations of pressure, different techniques have been proposed in the literature.
On the one hand, the \emph{supremizer} \citep{ballarin2015supremizer} and the Pressure Poisson Equation \citep{Stabile-SR-18}  introduce corrections either enriching the velocity space or performing a Helmholtz decomposition of the velocity, with the goal of obtaining a stable reduced system in the offline phase.
On the other hand, stability can be retrieved in the online phase via pressure-stabilised Petrov-Galerkin strategies \citep{Codina-BCI-14,Iliescu-CIVS-14}.
In other words, since the inf-sup condition is challenging to guarantee in the reduced-order framework, different alternatives have been devised to avoid the appearance of instabilities associated with the reduction of the approximation space for velocity and its consequent incompatibility with the reduced space for pressure.

In this work, a two-step procedure inspired by the consideration that stable approximations for mixed formulations can be obtained selecting a discrete space for velocity \emph{richer enough} with respect to the discrete space for pressure \citep{BREZZI199027,Bathe-CB-93} is followed.
First, the reduced basis for velocity is constructed setting two tolerances, $\varepsilon_u$ for $\bu$ and $\varepsilon_p$ for $p$, in the stopping criterion \eqref{eq:criterion}.
To guarantee that the subspace generated by the reduced basis for velocity is more accurate than the one for pressure, $\varepsilon_u$ is selected lower than $\varepsilon_p$. 
Indeed,
the reduced basis with low resolution selected for the pressure is sufficient to enforce the incompressible nature of the velocity but it does not provide an approximation of the pressure field with enough quality and physical significance. 
Hence, the second step of the procedure performs a simple post-process of the reduced velocity obtained with the POD-RB to recover a physically consistent pressure field.

More precisely,  note that the linear system to be solved at each iteration is a linearization of the nonlinear system
\begin{equation}\label{eq:NonLinearSystem}
\begin{bmatrix}
\bK(\nodalu) & \bm{B}^\TT \\
\bm{B} & \bm{0} 
\end{bmatrix}
\begin{bmatrix}
\nodalu  \\
\nodalp
\end{bmatrix} = 
\begin{bmatrix}
\nodalfu  \\
\bm{0} 
\end{bmatrix}.
\end{equation} 
Following the idea of velocity-pressure splitting \citep{Temam-84,GreshoSani1987}, the discretized momentum equation (i.e., the first row of \eqref{eq:NonLinearSystem})
\begin{equation}\label{eq:momentumNSsplit}
\bm{K}(\nodalu) \, \nodalu + \bm{B}^\top \nodalp = \nodalfu 
\end{equation}
can be used to compute the pressure corresponding to the velocity field $\nodalu$. Thus, pre-multiplying \eqref{eq:momentumNSsplit} by $\bm{B}$ yields to
\begin{equation}
\bm{B}  \bm{B}^\top \nodalp = \bm{B} (\nodalfu  - \bm{K}(\nodalu) \nodalu ),
\end{equation}
which is used as a simple post-process to recover the reduced pressure, once the reduced velocity is computed.
Alternative strategies for pressure recovery within incompressible Navier-Stokes are analyzed in  \cite{kean2020pressure}.

\section{Data augmentation}
\label{sec:data_augmentation}

As previously mentioned, a bottleneck of the standard POD-RB strategy is the necessity of having a sufficient number of snapshots, representative of the parametric family of solutions. Each snapshot is obtained solving a full-order problem, therefore it is demanding in terms of computational effort.  The idea of data augmentation consists of devising simple operations among the existing snapshots to create new functions that do likely approximate other snapshots, typically corresponding to intermediate parametric values.

In particular, in \cite{diez2021nonlinear}, a methodology was introduced to create new artificial snapshots using simple operations upon the original ones and a priori knowledge of the physical behaviour of the solution. 
For convection-diffusion problems, the operation was simply a pairwise product of the solution fields (a Hadamard, component by component, product of the vectors of nodal values). This is a sensible choice in that context because the product of pulses that propagate and diffuse tends to create intermediate pulses, which are exactly the missing snapshots. The strategy of augmenting the family of snapshots with these new elements was denoted as \emph{quadratic enrichment}.

In the present work, the data augmentation rationale is generalized to parametric problems in fluid mechanics involving incompressible flows. 
The two challenges of these problems (i.e., the incompressibility constraint and the nonlinearity of the convective field) are directly related to the underlying physics of the flow system, described by the conservation laws in \eqref{eq:NavierStokes}.
Hence, in order to introduce new, relevant information into the dataset with the artificial snapshots, it is of outmost importance that they represent the behavior of the fluid correctly.

The strategy to construct artificial snapshots features two steps. 
First, a list of possible pairwise combinations of the original snapshots is determined (Section \ref{sec_PairingSnapshots}).
Then, for each pairwise combination, different strategies are considered to design artificial snapshots fulfilling the balance laws by construction.
More precisely, the discussed approaches propose to enforce: 
\begin{itemize}
\item conservation of mass (Section \ref{sec_SolGeomAver});
\item conservation of mass and conservation of momentum (Section \ref{sec:Oseen}).
\end{itemize}

\subsection{Pairing of snapshots}
\label{sec_PairingSnapshots}

Data augmentation procedures rely on operations on previously computed snapshots.
In this work, a uniform sampling of the parametric space is considered to construct the initial dataset, without incorporating any problem-specific information. 

The augmentation procedure is based on pairwise combinations of solutions.
Note that in training sets with a significant number of snapshots (e.g., for high-dimensional problems), 
it can be computationally advantageous to 
combine only pairs of snapshots that are \emph{close} to one another,
by introducing an appropriate metric, instead of accouting for all possible combinations. 
More precisely,  in this work, the selection of pairs of snapshots to be combined is based on their proximity in the parametric space.
For the numerical experiments in Section \ref{sec:NumEx}, the Euclidean distance is employed but alternative solutions have been proposed in the literature, see, e.g., \cite{Aggarwal-AHK-01}.

\begin{figure}[]
	\centering{
		\begin{subfigure}{.35\textwidth}
			\includegraphics[width=\linewidth]{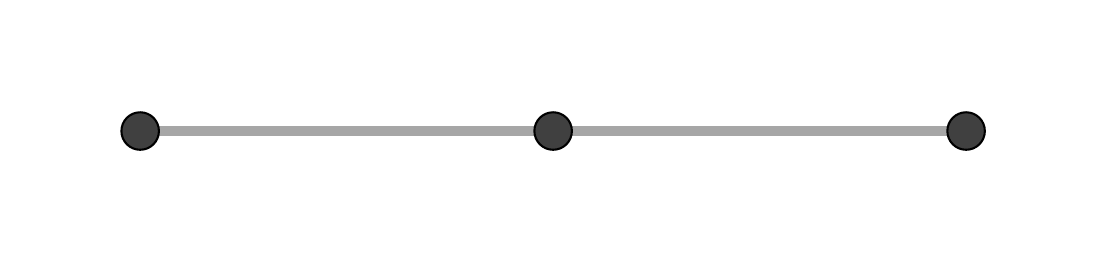}
			\caption{$\np =1$, $\nS = 3$.}
		\end{subfigure}\hspace{1cm}
		\begin{subfigure}{.35\textwidth}
			\includegraphics[width=\linewidth]{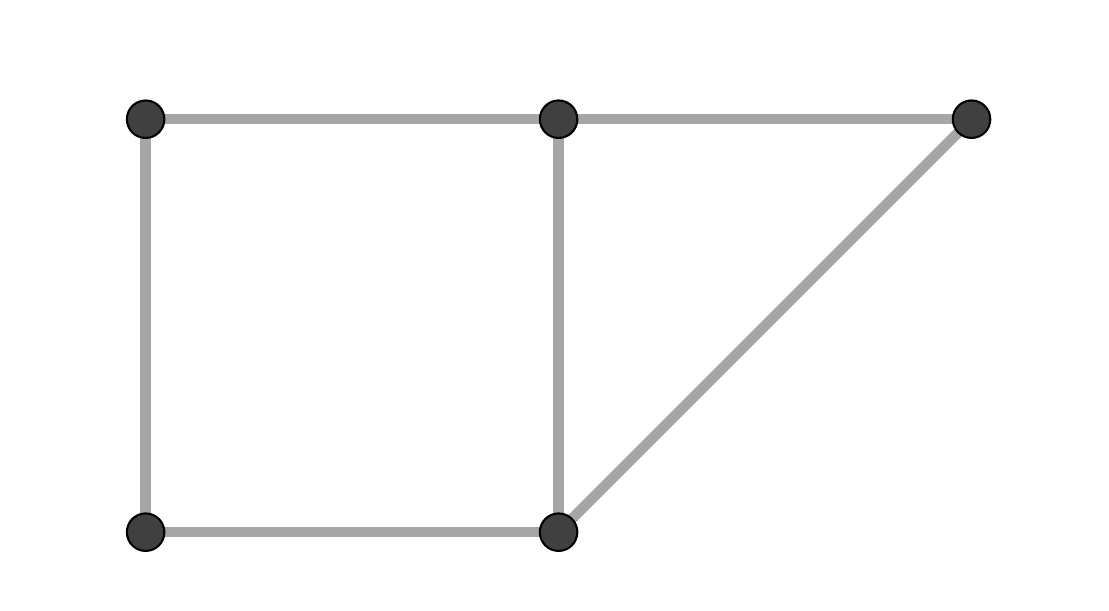}
			\caption{$\np =2$, $\nS = 5$.}
			\label{fig:pairing_2d}
		\end{subfigure}
		\caption{Pairing strategy of uniformly sampled snapshots in the parametric space.  
		(a) One parameter, three snapshots.
		(b) Two parameters, five snapshots.
		Each point corresponds to one snapshot in the initial dataset. 
		Lines denote the pairing of snapshots.  
		Each snapshot is paired with its nearest neighbor in each direction and the number of lines denote the number of combinations to be employed to generate artificial snapshots.}
		\label{fig:pairing}}
\end{figure}

Every snapshot of the dataset is thus paired with its nearest neighbor along each direction of the parametric space.
When two potential neighbors are equidistant, the algorithm selects the one that has not been previously paired to any other snapshot in order to maximize the information introduced during the design of the artificial data. 
In a parametric space of dimension $\np$, the algorithm thus identifies at least $\np$ pairs of snapshots per element of the initial training set. That is, the approach features a number of combinations scaling with the dimensionality of the problem.
The obtained pairings for a one-dimensional and a two-dimensional parametric space are exemplified in Figure \ref{fig:pairing}.

Finally, for each identified pair of snapshots, the procedures in Sections \ref{sec_SolGeomAver} and \ref{sec:Oseen} are executed to generate the artificial snapshots.
Alternative strategies, outside the scope of the present work, could be developed to further reduce the number of coupled snapshots during the augmentation procedure.

\subsection{Solenoidal field recovered from geometric average of stream functions}
\label{sec_SolGeomAver}

In this section, a strategy to expand the approximation space with meaningful new snapshots $\bu^\star$, produced with a cost-effective procedure, is presented. In order to be physically meaningful, the artificial snapshots are required to be incompressible, that is, to fulfill the mass conservation equation.

Consider two velocity snapshots $\bu_i$ and $\bu_j$, corresponding to parameters $\bmu_i$ and $\bmu_j$. The idea is to use them to easily compute a new velocity field $\bu^\star$. 
The \emph{quadratic enrichment} introduced in \cite{diez2021nonlinear} does not represent a viable option in this context as it assumes snapshots to have only positive values, whereas velocity fields can present both positive and negative components.  Moreover, the product of two snapshots would break the solenoidal character of the solution, no longer fulfilling the mass conservation equation.

To preserve incompressibility, a new snapshot is constructed starting from the information of the stream function. First, the velocity fields $\bu_i$ and $\bu_j$ are expressed in terms of the corresponding stream functions $\Psi_i$ and $\Psi_j$.  The stream functions $\Psi_i$ and $\Psi_j$ are thus geometrically averaged by computing a new stream function
\begin{equation}\label{eq:solenoidal_av}
\Psi^\star = (\Psi_i) ^\alpha   (\Psi_j) ^{1-\alpha} ,
\end{equation}
for a scalar weight $\alpha \in ]0,1[$.
Note that, in order to guarantee that $\Psi^\star$ is well defined, $\Psi_i$ and $\Psi_j$ are appropriately translated to ensure that their corresponding values are positive.
Then, the velocity field $\bu^\star$ corresponding to $\Psi^\star$ is an inexpensive, artificial snapshot, expected to be an approximation of the solution corresponding to some intermediate parameter, between $\bmu_i$ and $\bmu_j$. 

Note that the operations of recovering $\Psi_i$ and $\Psi_j$ from $\bu_i$ and $\bu_j$ and then computing $\bu^\star$ from $\Psi^\star$ are standard and computationally affordable, see \cite[Chapter~6]{Kundu-book}. Indeed, in 2D, $\Psi_i$ is computed from $\bu_i$ by solving the Poisson problem
\begin{equation}\label{eq:SolPoisson}
\left\{
\begin{aligned}
-\grad^2 \Psi_i &= \frac{\partial [\bu_i]_y}{\partial x} - \frac{\partial [\bu_i]_x}{\partial y}  & \text{ in } \Omega,  \\
\grad \Psi_i \cdot \bn &= (-[\bu_i]_y,[\bu_i]_x)\cdot\bn  & \text{ on } \partial \Omega,
\end{aligned}
\right.
\end{equation}
where $[\bu_i]_x$ and $[\bu_i]_y$ stand for the first and second components of $\bu_i$, respectively.
Note that only Neumann boundary conditions are set. Therefore, to have a unique solution, the value of $\Psi_i$ is to be set at some arbitrary point (e.g., $\Psi_i=0$ at a node lying on the boundary $\partial \Omega$). The velocity field $\bu^\star$ corresponding to the new stream function $\Psi^\star$ is thus readily computed, see e.g.  \cite[Chapter~6]{Kundu-book}, as
	\begin{equation}\label{eq:VelFromPsi}
	\bu^\star 
	= 
	\begin{bmatrix}
	\displaystyle\frac{\partial \Psi^\star}{\partial y} \\[1.1ex]
	-  \displaystyle\frac{\partial \Psi^\star}{\partial x}
	\end{bmatrix}
	= 
	\begin{bmatrix}
	\,\,\, 0 & 1\\
	-1 & 0
	\end{bmatrix}
	\grad \Psi^\star 
	.
	\end{equation}
Taking different values of $\alpha$ between 0 and 1, different intermediate velocity fields are generated from each pair of actual snapshots in the original training set. 
It is worth noticing that, regardless of the value of $\alpha$ selected for the construction of the function $\Psi^\star$, the corresponding velocity $\bu^\star$ obtained from \eqref{eq:VelFromPsi} is solenoidal by construction because
$$
\grad \cdot \bu^\star = \frac{\partial^{2} \Psi^\star}{\partial x \partial y} - \frac{\partial^{2} \Psi^\star}{\partial y \partial x}=0 ,
$$
assuming that $\Psi^\star$ is sufficiently regular.

In order to ensure incompressibility, artificial velocity snapshots $\bu^\star$ are thus generated starting from the intermediate stream functions $\Psi^\star$.
The strategy to create the intermediate stream functions $\Psi^\star$ introduced in \eqref{eq:solenoidal_av} is the simplest one producing linearly independent snapshots but alternative formulations guaranteeing that $\Psi^\star$ is linearly independent of $\Psi_{i}$ and $\Psi_{j}$ might also be admissible. 
Indeed,  note that if the different $\Psi^\star$ were linearly dependent on $\Psi_{i}$ and $\Psi_{j}$, then $\bu^\star$ would be linearly dependent on $\bu_{i}$ and $\bu_{j}$, and the expanded basis would only introduce redundant information.

As it is shown in the examples, in some cases the new snapshots generated with this technique are not bringing into the training set new pieces of information that contribute to improving the reduced-order solution. 
Thus, in order to further enrich the quality of the basis, a correction imposing both mass conservation and momentum balance via the Oseen equation is proposed.

\subsection{Physics-informed enhancement using Oseen equation}
\label{sec:Oseen}

The data augmentation strategy discussed in this section aims to concurrently enforce mass and momentum conservation. To achieve this goal with a cost lower than the computation of a new snapshot, the proposed method relies on solving a linearized version of equation \eqref{eq:NavierStokes}. 
More precisely, let $\bu^\star$ be an artificial snapshot created following the strategy presented in Section \ref{sec_SolGeomAver}. By construction, it fulfills the mass conservation equation.
The resulting velocity $\bu^\star$ is thus employed as a \emph{known} convective field in the momentum conservation equation for a new, artificial snapshot $\bu^{\star\star}$ enchancing the compliance with the physics of the flow system via the fulfilment of the balance of forces.

In particular, the solenoidal field  $\bu^\star$ is used as input data for an Oseen problem, which is a linearization of the steady-state incompressible Navier-Stokes equation reading
\begin{equation}\label{eq:Oseen}
\left\{
\begin{aligned}
-\grad \cdot \left( \nu \grad \bu^{\star\star} - p \bI \right) +  \left( \bu^\star \cdot \grad \right) \bu^{\star\star} &= \bm{0}  &\text{in } \Omega, \\
\grad \cdot \bu^{\star\star} &= 0  & \text{in } \Omega,
\end{aligned}
\right.
\end{equation}
complemented with the same type of boundary conditions as \eqref{eq:NavierStokes}.
In the case of boundary conditions which depend on the parameters of the problem $\bmu$, an approximation of the corresponding value is to be inferred from $\bu^\star$. 

Note that the discrete form of \eqref{eq:Oseen} is a variation of \eqref{eq:NonLinearSystem}, namely
\begin{equation}\label{eq:discreteOseen}
\begin{bmatrix}
\bK(\nodalu^\star) & \bm{B}^\TT \\
\bm{B} & \bm{0} 
\end{bmatrix}
\begin{bmatrix}
\nodalu^{\star\star}  \\
\nodalp
\end{bmatrix} = 
\begin{bmatrix}
\nodalfu  \\
\bm{0} 
\end{bmatrix} ,
\end{equation}
with the peculiarity that the resulting problem \eqref{eq:discreteOseen} is linear. 
Of course,  to construct the Oseen enhancement enforcing the momentum balance,  this problem needs to be solved multiple times. Indeed,  for each artificial snapshot $\bu^\star$ created as solenoidal average of any pair of original snapshots, one problem \eqref{eq:discreteOseen} is to be solved. 
Although this might seem computationally demanding, only few terms in the matrix require to be updated when a new artificial snapshot is considered. In addition, exploiting the linearity of the resulting problem, efficient solution strategies can be devised using suitable Krylov subspace methods. A more detailed discussion on the efficiency of the proposed approach is presented in Section \ref{sec:Discussion}.

Following a similar rationale, an alternative physics-informed enrichment can be devised by substituting the solenoidal average in \eqref{eq:Oseen} with an artificial velocity field $\bu^\star$ obtained as a linear combination of two original snapshots. More precisely, given two snapshots $\bu_i$ and $\bu_j$,  an artificial velocity is first defined as $\bu^\star = \alpha \bu_i + (1 -\alpha) \bu_j$, with $\alpha \in ]0,1[$ being a weighting coefficient. Then, solving the Oseen equation \eqref{eq:Oseen}, a new artifical snapshot $\bu^{\star\star}$ is determined.
Numerical examples comparing the two approaches are discussed in Section \ref{sec:NumEx}.

This alternative generates a new family of linearly independent, artifical snapshots $\bu^{\star\star}$. Although the artificial velocity $\bu^\star$ is linearly dependent on the original ones, the artifical snapshot $\bu^{\star\star}$ obtained from the Oseen problem is independent and actually enriches the basis. 
Moreover,  the physics-based enrichment relying on an artificial velocity obtained as linear combination of two original snapshots allows a straightforward extension to the three-dimensional case,  where the notion of stream function is only available in the case of axisymmetric flows \cite[Chapter~6]{Kundu-book}.

\section{Numerical examples}
\label{sec:NumEx}

In this Section, the presented data augmentation strategies are tested to construct artificial snapshots within the POD-RB formulation, for the steady-state incompressible Navier-Stokes equations. Three benchmarks involving parametrized flow problems, in two and three dimensions, are considered, with one and two user-controlled parameters.

Full-order simulations are performed using the open-source library FEniCS \citep{FenicsProject,LangtangenLogg2017}. In all computations, the relative tolerance for the convergence of the Newton-Raphson scheme is set to $10^{-8}$ for both the pressure and the velocity fields. 

In the 2D examples, the three strategies discussed in Section \ref{sec:data_augmentation} are tested to enrich the velocity training set, by incorporating artificial snapshots computed by:
\begin{enumerate}[label=(\roman*)]
	\item solenoidal average of two snapshots (i.e., imposing mass conservation);
	\item solenoidal average of two snapshots, enhanced with the Oseen equation (i.e., imposing mass and momentum conservation);
	\item linear combination of two snapshots enhanced with the Oseen equation (i.e., imposing mass and momentum conservation).
\end{enumerate}
As mentioned in Section~\ref{sec:data_augmentation},  the enrichment procedure based on the solenoidal average of the snapshots relies on the computation of the stream function which is only suitable for two-dimensional or three-dimensional, axisymmetric flows.  Hence, only strategy (iii) is presented for the 3D test case in Section~\ref{sc:microswimmer}, which features a parametric Navier-Stokes flow in a space of dimension $5$.
It is worth noticing that enrichment strategy (iii) is applicable to parametric incompressible flows, independently of the number of spatial dimensions $d$ and the number of parameters $\np$.
Whilst the computational benefit of performing data augmentation is expected to be even more relevant for problems with $d+\np >5$,  these entail a series of technical developments, including the high-dimensional sampling \cite{Vono-VDC-22} to construct the initial training set and the accurate approximation of high-dimensional functions \cite{Brugiapaglia-ABW-22},  which lie beyond the scope of the present contribution.

\subsection{Jet control for flow past a cylinder}

\begin{figure}[h!]
	\centering
	\includegraphics[width=0.45\textwidth]{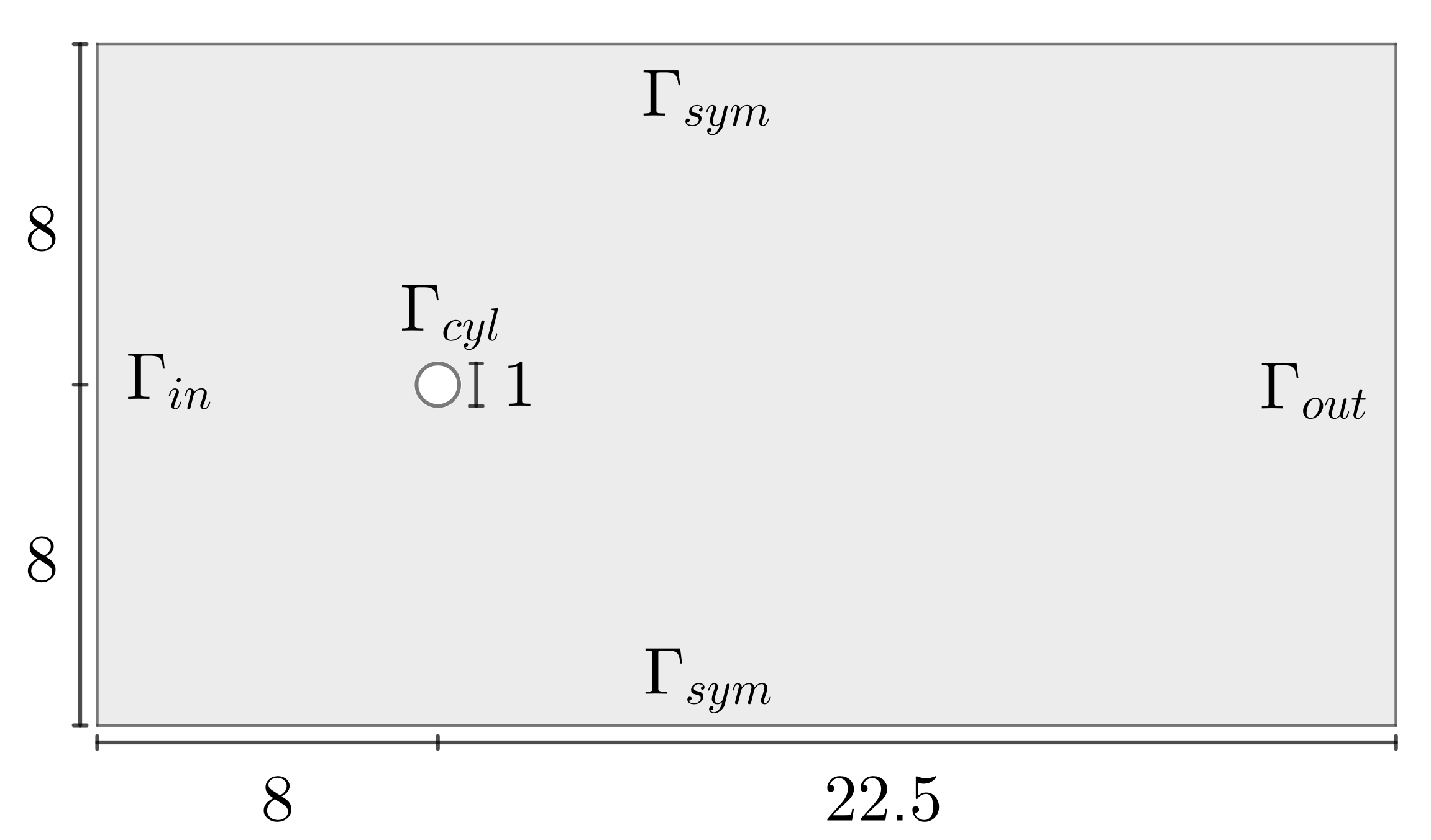}
	\hspace{5mm}
	\raisebox{2mm}{\includegraphics[width=0.3\textwidth]{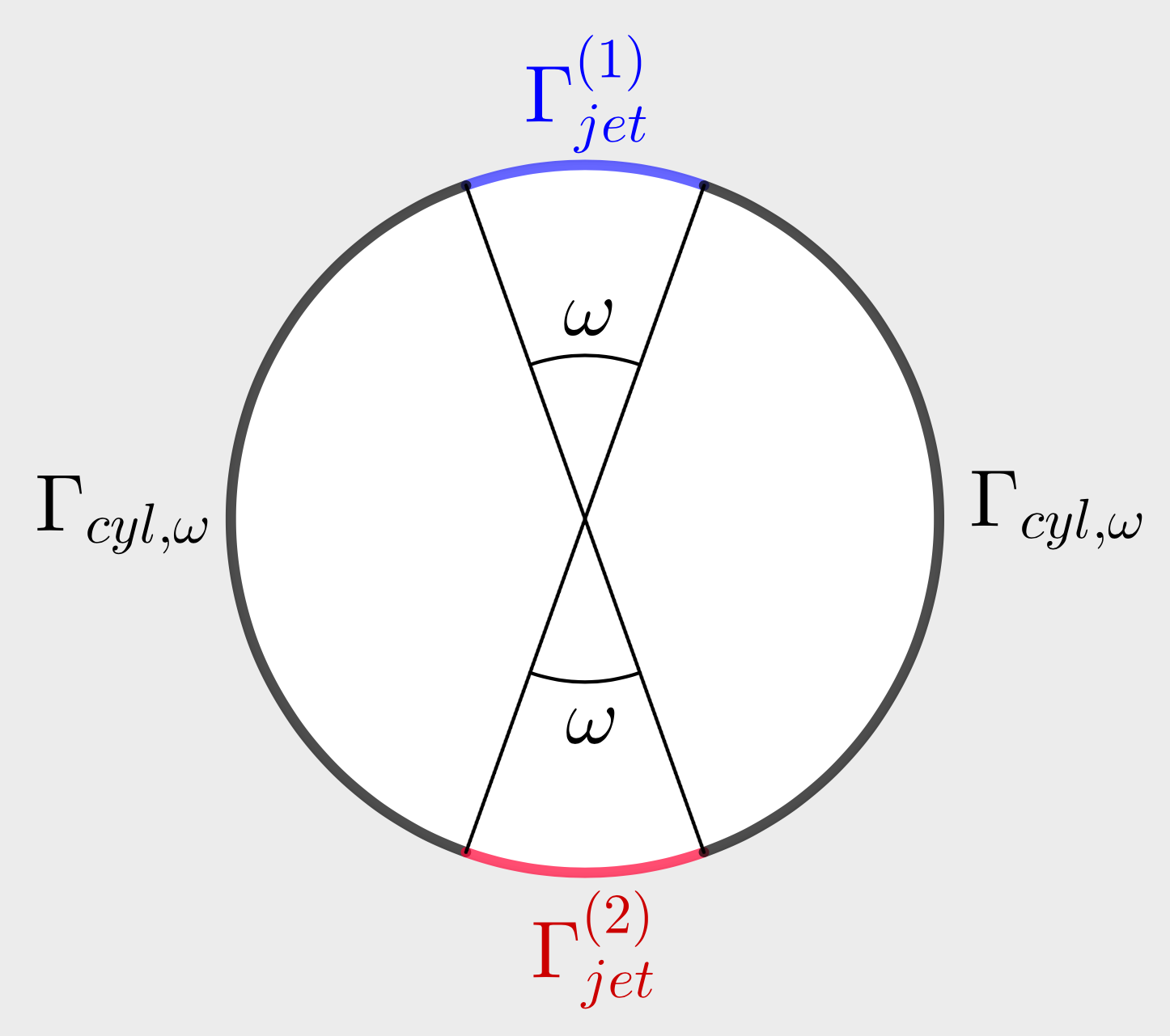}}
	\caption{\normalsize Geometry of the domain and boundary partition.}
	\label{fig:jets_domain}
\end{figure}

The first example, inspired by \cite{Rabault-RKJRC-19}, considers the flow past a cylinder where two blowing/sucking jets are introduced.
Consider the domain $\Omega = [0,30.5]\times[0,16] \setminus \mathcal{B}_{0.5}(8,8)$, where $\mathcal{B}_{\R}(x_c,y_c)$ denotes the two-dimensional cylinder of radius $\R$ centered in $(x_c,y_c)$. The boundary $\partial\Omega$ is partitioned according to the schematics reported in Figure \ref{fig:jets_domain}.

The fluid enters the domain with uniform horizontal velocity $\Uin$ through the inlet on the left boundary,  symmetry conditions are imposed on the top and bottom, and an outlet with homogeneous Neumann condition is added on the right vertical contour.
On the surface of the cylinder, no-slip conditions are enforced, except for two portions, on the upper and lower parts, where two jets are inserted (Figure \ref{fig:jets_domain}, right).
The boundary conditions read
\begin{subequations}
	\begin{equation}
	\begin{aligned}
	\bu &= [ \Uin, 0]  && \text{ on } \Gamma_{in}, \\
	\bu \cdot \bn &= 0  && \text{ on } \Gamma_{sym}, \\
	\bm{\tau} \cdot \left( \nu \grad \bu - p \bI \right)\bn &= 0 && \text{ on } \Gamma_{sym}, \\
	\bu &= \bm{0} && \text{ on } \Gamma_{cyl,w}, \\
	\bu &=  \gamma\UvecJet  && \text{ on } \Gamma_{jet}^{(1)}, \\
	\bu &= -\gamma\UvecJet && \text{ on } \Gamma_{jet}^{(2)}, \\
	\left( \nu \grad \bu - p \bI \right)\bn &= \bm{0}  && \text{ on } \Gamma_{out},
	\end{aligned}
	\end{equation}
	where $\bm{\tau}$ denotes the unit tangent vector,  $\gamma$ is a dimensionless parameter that modulates the maximum velocity of the jets, whereas their velocity profiles are given by
	\begin{equation}
	\UvecJet(x,y) = 2\Uin \cos \left( \frac{\pi}{\omega} \left(\theta - \frac{\pi}{2}\right)\right) \left[ x-8, y -8  \right],
	\end{equation}
\end{subequations}
with $\omega = 5\pi/36$ being the opening angle of each jet. 

The kinematic viscosity is set to $\nu = 0.01$ and the Reynolds number is $\Re = \Uin \D/\nu$, the characteristic length being the diameter $\D=1$ of the cylinder.
The full-order solver relies on a $\mathbb{P}_2/\mathbb{P}_1$ finite element discretization computed on a mesh with $10,962$ triangular elements, suitably refined in the boundary layer region and in the wake of the cylinder. The resulting discretization features a total of $44,400$ degrees of freedom for velocity and $5,610$ for pressure. 

\begin{figure}[h!]
	\centering{
	\begin{subfigure}{.4\textwidth}
		\includegraphics[width=\linewidth]{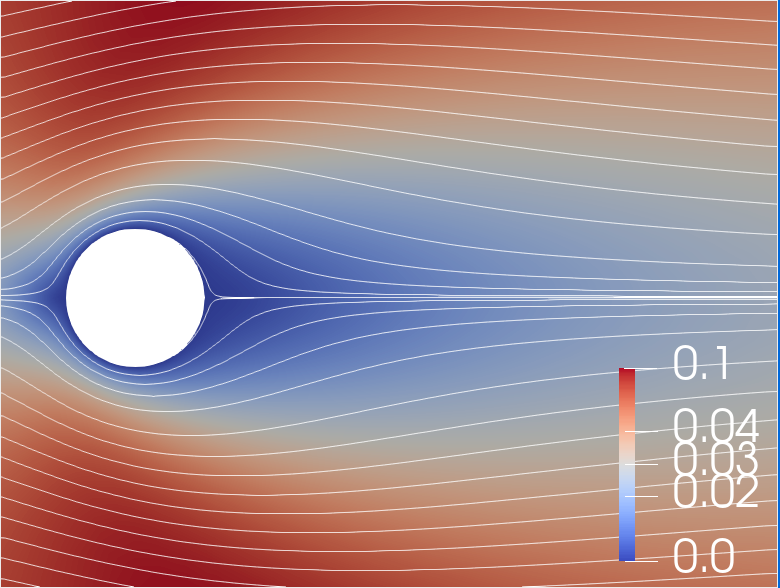}
		\caption{$\Re = 5$,  no jets ($\gamma = 0$).}
	\end{subfigure}\hspace{3mm}
	\begin{subfigure}{.4\textwidth}
		\includegraphics[width=\linewidth]{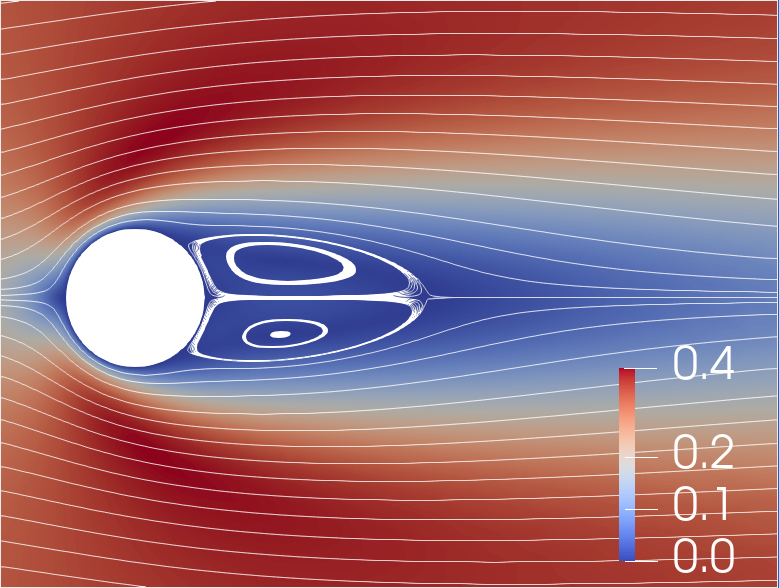}
		\caption{$\Re = 30$,  no jets ($\gamma = 0$).}
	\end{subfigure}
	\vspace{1mm}

	\begin{subfigure}{.4\textwidth}
		\includegraphics[width=\linewidth]{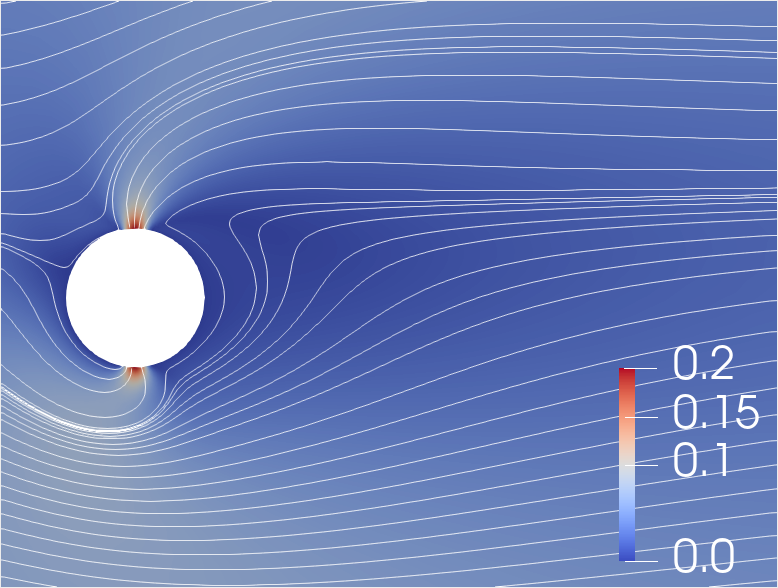}
		\caption{$\Re = 5$,  maximum jets ($\gamma = 4$).}
	\end{subfigure}\hspace{3mm}
	\begin{subfigure}{.4\textwidth}
		\includegraphics[width=\linewidth]{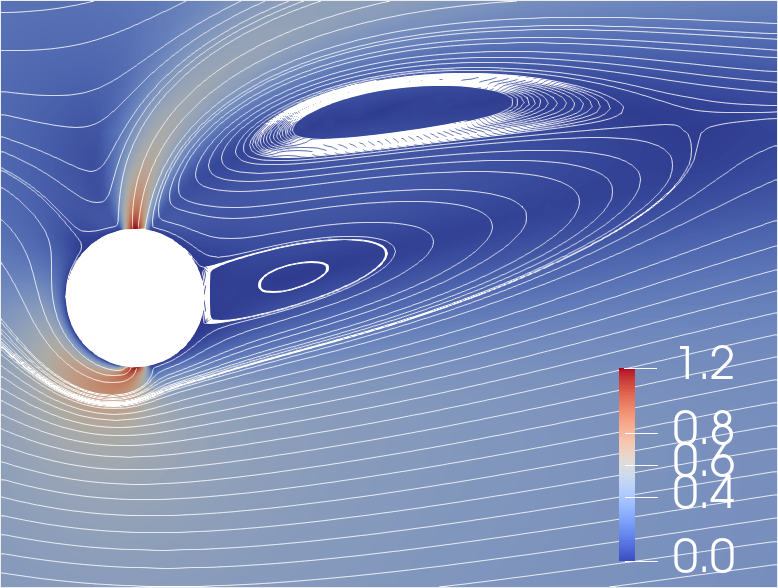}
		\caption{$\Re = 30$  maximum jets ($\gamma = 4$).}
	\end{subfigure}
	\vspace{1mm}

	\caption{Velocity streamlines in the wake of the cylinder (area of interest $[7,13]\times [6,10]$) for representative values of the parameters $\Re \in [5,30]$ and $\gamma\in[0,4]$ (jet velocity). The color scales indicate the velocity module.}
	\label{fig:jets_family}}
\end{figure}

The problem under analysis features two parameters: the Reynolds number $\Re \in [5,30]$ (through the inlet velocity $\Uin$) and the maximum velocity of the jets controlled by $\gamma \in [0,4]$.  The velocity field for a representative set of solutions is shown in Figure \ref{fig:jets_family}.

Truncation tolerances in POD are set to $\varepsilon_u = 10^{-3}$ for the velocity training set, and $\varepsilon_p = 0.25$ for the pressure training set. Recall that the pressure is later reconstructed as a post-process of the POD-RB velocity, as described in Section \ref{sec:Stability}.

Data augmentation is performed with weighting coefficients $\alpha = 0.1, 0.2, \dots, 0.9$ to construct $\bu^\star$, both for the solenoidal averages, (i) and (ii), and for the linear combination, (iii), strategies. This leads to the generation of $9$ artificial snapshots per each pair of originally-computed snapshots.

\subsubsection{Single-parameter analysis}

First, the proposed data augmentation strategies are analyzed in the case of a single parameter. Considering the parametrized Reynolds number $\Re$, the value $\gamma = 4$ is fixed. 
The original training set for the POD consists of only two snapshots, associated with $\Re = 5$ and $\Re = 30$. 
Upon augmentation, an enriched velocity training set of $11$ snapshots is obtained.

\begin{figure}[h!]
	\centering
	\includegraphics[width=0.6\textwidth]{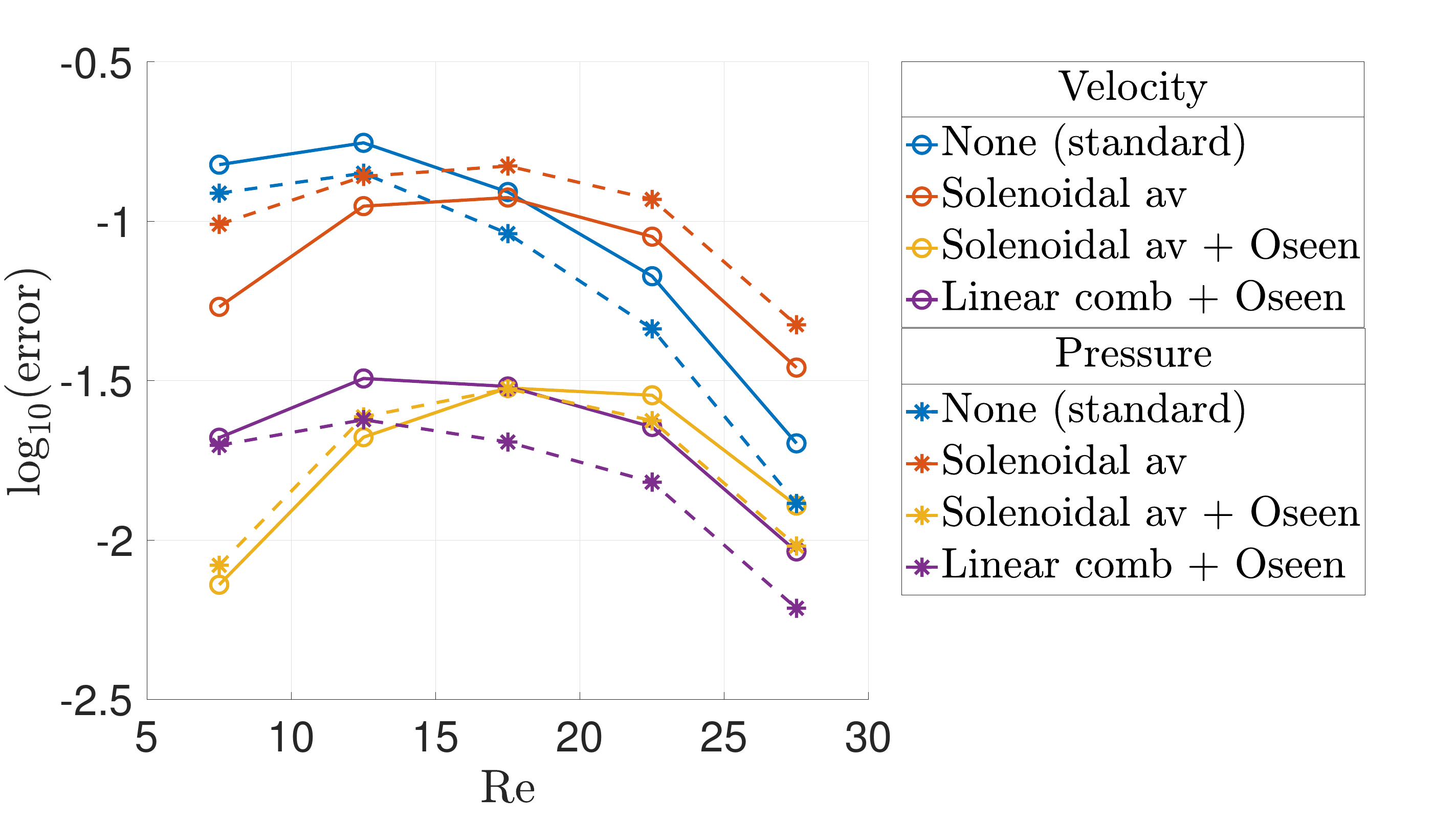}
	
	\includegraphics[width=0.6\textwidth]{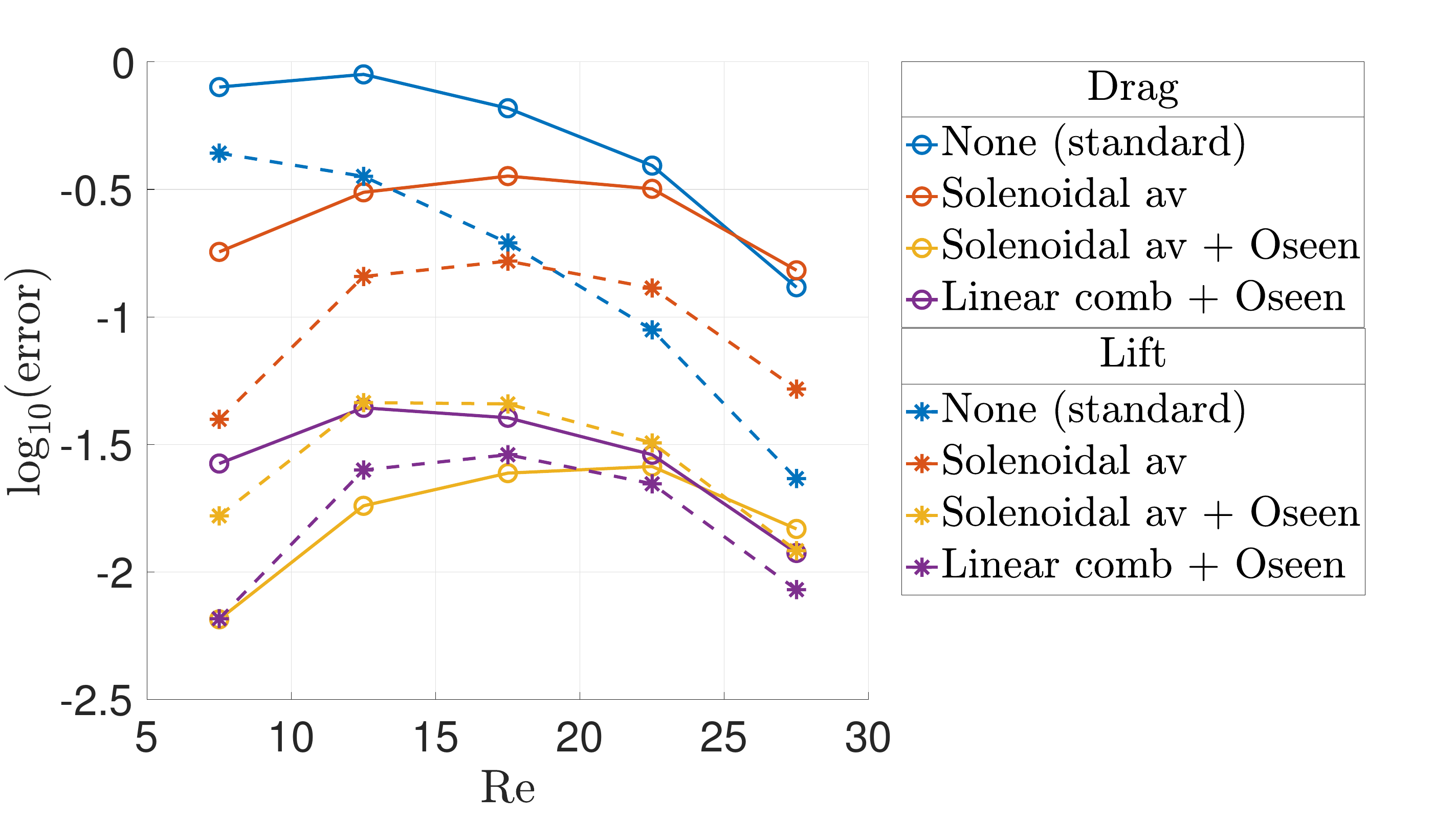}
	
	\caption{Relative errors of the POD-RB approximation, measured in the Euclidean norm, of velocity and pressure (top), drag and lift (bottom) for different strategies of data augmentation, with parametrized $\Re \in [5,30]$ and fixed $\gamma = 4$. }
	\label{fig:jets_errors_Re}
\end{figure}

The reduced dimension for the pressure is $\nkp = 1$. For the velocity, the reduced dimensions are $\nku = 1$ for the standard POD, 
$\nku = 3$ in the case of solenoidal averages enrichment, and $\nku = 4$ in the two cases with Oseen enhancement. The relative errors, measured in the Euclidean norm, of the reduced-order approximation at intermediate values of the Reynolds number are reported in Figure \ref{fig:jets_errors_Re}. 
In general, the results show that enriching the training set exclusively by means of solenoidal averages, does not improve the performance of the POD-RB.  In particular, the linear combination of the original snapshots inherent to the POD seems to contain the same information provided by the new, artificial snapshots obtained enforcing the mass conservation principle.
On the contrary, the physics-informed Oseen enhancement, in which both mass and momentum conservation are enforced, leads to more accurate solutions: errors improve for velocity, pressure, drag and lift in all tested values. For the velocity and pressure, the improvement is more significant for lower values of $\Re$, around one order of magnitude, and is less pronounced as $\Re$ increases and the associated relative errors decrease. The same tendency on error evolution is observed for the lift. Regarding the drag, errors improve approximately one order of magnitude for all values of $\Re$. 

\begin{figure}[h!]
	\centering
	\includegraphics[width=0.6\textwidth]{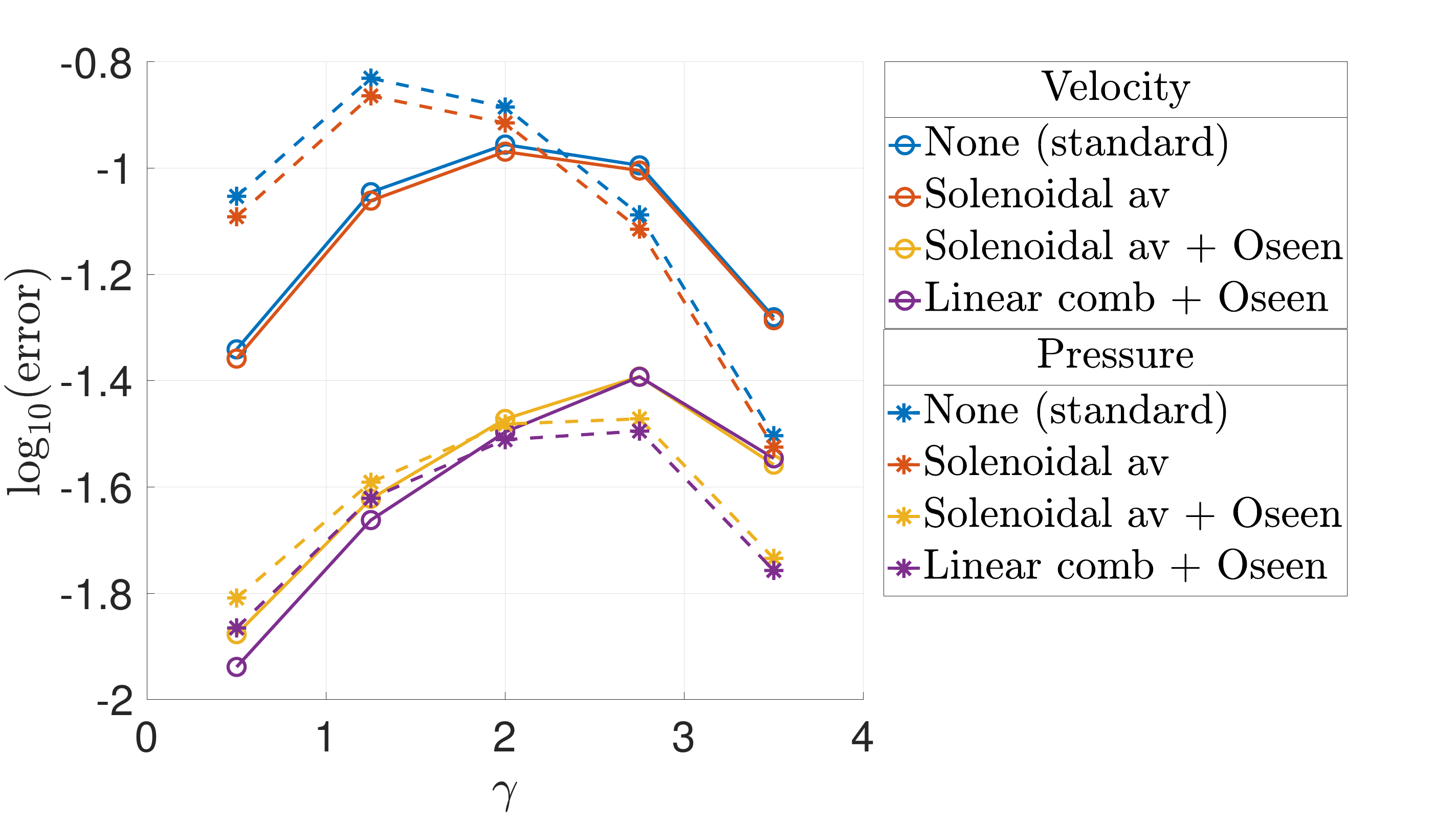}\\	
	\includegraphics[width=0.6\textwidth]{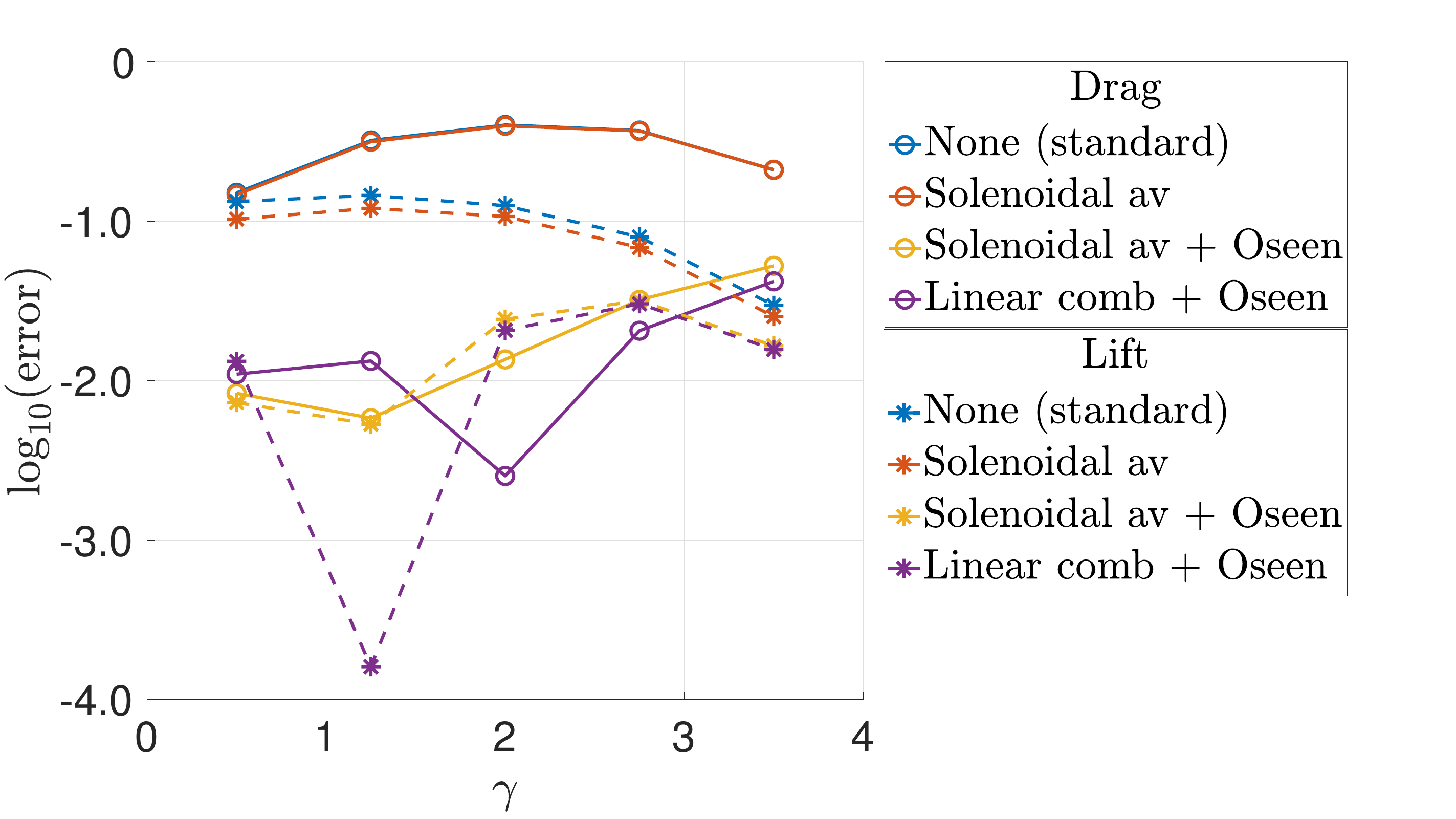}
	\caption{Relative errors of the POD-RB approximation, measured in the Euclidean norm,  of velocity and pressure (top), drag and lift (bottom) for different strategies of data augmentation, with parametrized $\gamma \in [0,4]$ and fixed $\Re = 30$ . 
	}
	\label{fig:jets_errors_alpha}
\end{figure}

The same test is reproduced fixing $\Re = 30$ and considering $\gamma$ as unique parameter.  As before, the training set consists of two snapshots, corresponding to $\gamma = 0$ and $\gamma = 4$, whereas the enriched velocity training set contains $11$ snapshots.
It follows that the reduced dimension for pressure is $\nkp = 1$, whereas for velocity the standard POD-RB approximation features a reduced dimension $\nku = 1$ which increases to $\nku = 4$ for the case of enrichment with solenoidal averages and $\nku = 6$ for the  cases with Oseen enhancement. 
Figure \ref{fig:jets_errors_alpha} shows the relative errors at some intermediate values of $\gamma$, measured in the Euclidean norm. 
On the one hand, negligible differences are observed in the POD-RB approximation between the standard ROM and the solenoidal averages enrichment strategy, highlighting that only imposing mass conservation during the construction of the artificial snapshots does not introduce relevant information not present in the original dataset.
On the other hand, the Oseen enhancement is able to incorporate additional information, stemming from momentum conservation, into the approximation space and leads to more accurate results.

\begin{remark}
The two data augmentation approaches exploiting both mass and momentum conservation principles lead to similar results for the POD-RB in this example.
Indeed,  since the range of the Reynolds number under analysis is fairly low,  all flows are highly viscous and vorticity effects are limited, whence both solenoidal averages and linear combinations of snapshots are able to provide sufficient information on the balance of forces.
The performance of these two strategies in the context of transient problems and flows with more relevant vorticity effects should be thoroughly examined in future works.
\end{remark}

\subsubsection{Two-parameter analysis}

Finally, the concurrent parametrization of $\Re$ and $\gamma$ is studied, with a training set comprising $30$ snapshots obtained from the uniform discretization of the parametric space. 
In particular, we consider the local POD-RB approximation of the problem, accounting for the $5$ closest snapshots in the parametric space. 
Upon data augmentation, new artificial velocities are included in the training set. 
Starting from a uniform discretization of $5$ snapshots, $6$ pairs of snapshots are determined, as sketched in Figure \ref{fig:pairing_2d}. Following the strategy described in the previous subsection, $9$ artificial snapshots are thus computed per each pair corresponding to the values of $\alpha = 0.1, 0.2, \dots, 0.9$ as interpolation weight. The process results in a collection of a total of $59$ snapshots, $5$ from the original dataset and $54$ artificial ones. 

In Table \ref{tab:2param},  a comparison of the different enrichment techniques is presented. The dimension of the POD reduced spaces after data augmentation is reported, highlighting the capability of the proposed approaches to generate artificial snapshots.
For a new point, corresponding to the parameters $\Re^\star = 17.5$ and $\gamma^\star = 1.5$,  not included in the training set, the relative errors of the reduced-order approximations, measured in the Euclidean norm, are also listed. 
The same behavior observed in the single-parameter analysis is confirmed in the case of two parameters. Data augmentation only enforcing mass conservation via solenoidal averages leads to similar approximation errors as the standard POD. On the contrary, accounting for both mass and momentum conservation by means of the Oseen enhancement strategy allows to significantly improve the errors for all the measured variables.
\begin{table}[h!]
	\centering
	\begin{tabular}{ lcccccc } 
		\hline
		Data  & \multicolumn{2}{c}{Dimension} & \multicolumn{4}{c}{Error} \\[-0.5em]
		augmentation & $\nku$ & $\nkp$ & Velocity & Pressure  & Drag & Lift \\
		\hline 
		None & $4$ & $2$ & $ 3.58\cdot 10^{-3}$ & $ 6.48\cdot 10^{-3}$ & $ 9.39\cdot 10^{-3}$ & $ 2.06\cdot 10^{-3}$  \\ 
		(i) & $7$ & $2$ & $3.66\cdot 10^{-3}$ & $6.36\cdot 10^{-3}$ & $  8.99\cdot 10^{-3}$ & $4.63\cdot 10^{-3}$   \\ 
		(ii)  & $8$ & $2$ & $5.97\cdot 10^{-4}$ & $1.35\cdot 10^{-3}$ & $3.07\cdot 10^{-4}$ & $1.22\cdot 10^{-3}$   \\ 
		(iii) & $7$ & $2$ & $9.80\cdot 10^{-4}$ & $1.60\cdot 10^{-3}$ & $6.42\cdot 10^{-4}$ & $1.92\cdot 10^{-3}$  \\ 
		\hline
	\end{tabular}
	\caption{Relative errors of the POD-RB approximation, measured in the Euclidean norm,  at $\Re^\star = 17.5$ and $\gamma^\star = 1.5$ for different strategies of data augmentation: none, standard POD; (i) solenoidal average; (ii) solenoidal average with Oseen enhancement; (iii) linear combination with Oseen enhancement. 
	}
	\label{tab:2param}
\end{table}

\subsection{Lid-driven cavity with parametrized jets}

\begin{figure}[h!]
	\centering
	\includegraphics[width=0.3\textwidth]{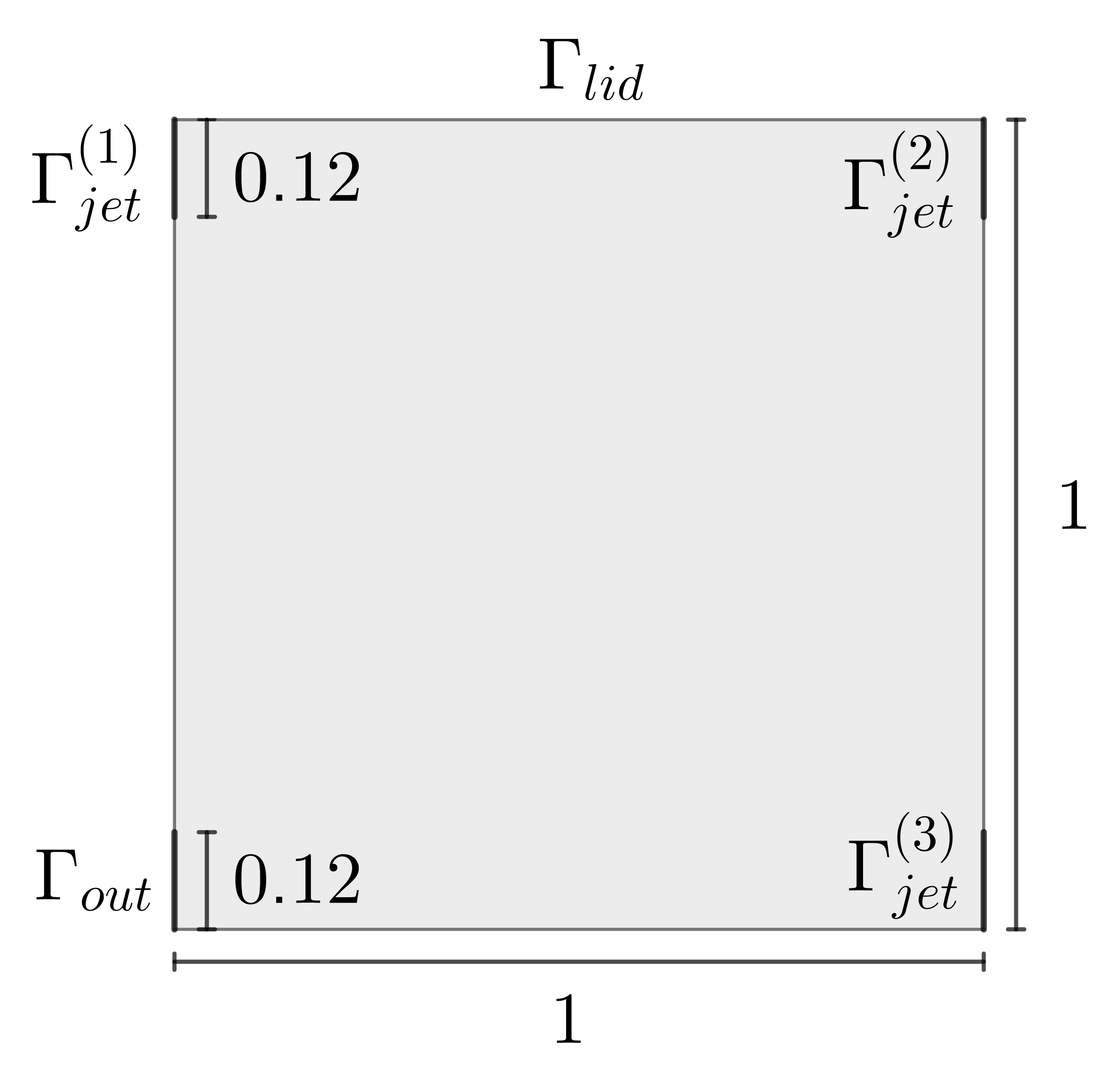}
	
	\caption{\normalsize Geometry of the domain and boundary partition.}
	\label{fig:lid_setting}
\end{figure}

In this Section,  the test case of the lid-driven cavity with parametrized jets in \cite{TsiolakisPGD2020} is considered. In this problem, three jets are introduced into the classical non-leaky lid-driven cavity \citep{ghia1982highRe}.
The spatial domain $\Omega = [0,1]^2$ is depicted in Figure \ref{fig:lid_setting}.  

A horizontal velocity profile $\Ulid(x)$ is imposed on the top lid of the domain, $\Gamma_{lid}$, and the three jets with horizontal velocity profiles $\mu \Ujet^{(i)}(y)$ are introduced on the vertical walls at the locations defined by
\begin{equation}
\begin{aligned}
\Gamma_{jet}^{(1)} &= \{ x = 0,  \ y\in[0.88,1] \} \, , \\
\Gamma_{jet}^{(2)} &= \{ x = 1,  \ y\in[0.88,1] \} \, , \\
\Gamma_{jet}^{(3)} &=  \{ x = 1,  \ y\in[0,0.12] \} \, .
\end{aligned}
\end{equation} 
Moreover,  an outlet surface is added on $\Gamma_{out} = \{  x = 0,  \ y\in[0,0.12] \}$.

The complete set of imposed boundary conditions read
\begin{subequations}
	\begin{equation}
	\begin{aligned}
	\bu &= [\Ulid,0] && \textnormal{ on } \Gamma_{lid} , \\
	\bu &= [\mu \Ujet^{(i)},0] && \textnormal{ on } \Gamma_{jet}^{(i)},  \ i=1,\ldots,3 ,\\
	\left( \nu \grad \bu - p \bI \right)\bn &= \bm{0} && \textnormal{ on } \Gamma_{out}, \\
	\bu &= \bm{0} && \textnormal{ elsewhere} ,
	\end{aligned}
	\end{equation}
	where the horizontal component of the lid velocity is given by
	\begin{equation}
	\Ulid(x) = \begin{cases}
	10 x/0.06 & \textnormal{for } x\in[0,0.06], \\
	10 & \textnormal{for } x\in(0.06,0.94), \\
	10(1-x)/0.06 & \textnormal{for } x\in[0.94,1] ,
	\end{cases}
	\end{equation}
	whereas the profiles $\Ujet^{(i)}$ on $\Gamma_{jet}^{(i)},  \ i=1,\ldots,3$ of the horizontal components of the jet velocities are defined as
	\begin{equation}
	\begin{aligned}
	\Ujet^{(1)}(y) &= 1-\cos\left(2\pi (y-0.88)/0.12\right) , \\
	\Ujet^{(2)}(y) &= 1-\cos\left(2\pi (y-0.88)/0.12\right), \\
	\Ujet^{(3)}(y) &= -1+\cos\left(2\pi y/0.12\right) , 
	\end{aligned}
	\end{equation}
	and scaled by the parameter $\mu\in(0,1]$.
\end{subequations}

The kinematic viscosity is set to $\nu = 0.01$ and the resulting Reynolds number computed using $\max | \Ulid |$ as characteristic velocity is $\Re = 1,000$. 
The high-fidelity solutions are computed using $\mathbb{P}_2/\mathbb{P}_1$ finite element basis on a mesh of $20,000$ triangular elements. The resulting discretization consists of $80,802$ degrees of freedom for velocity and $10,201$ for pressure. 

\begin{figure}[h!]
	\centering
	\begin{subfigure}[b]{0.305\textwidth}
		\centering
		\includegraphics[width=\textwidth]{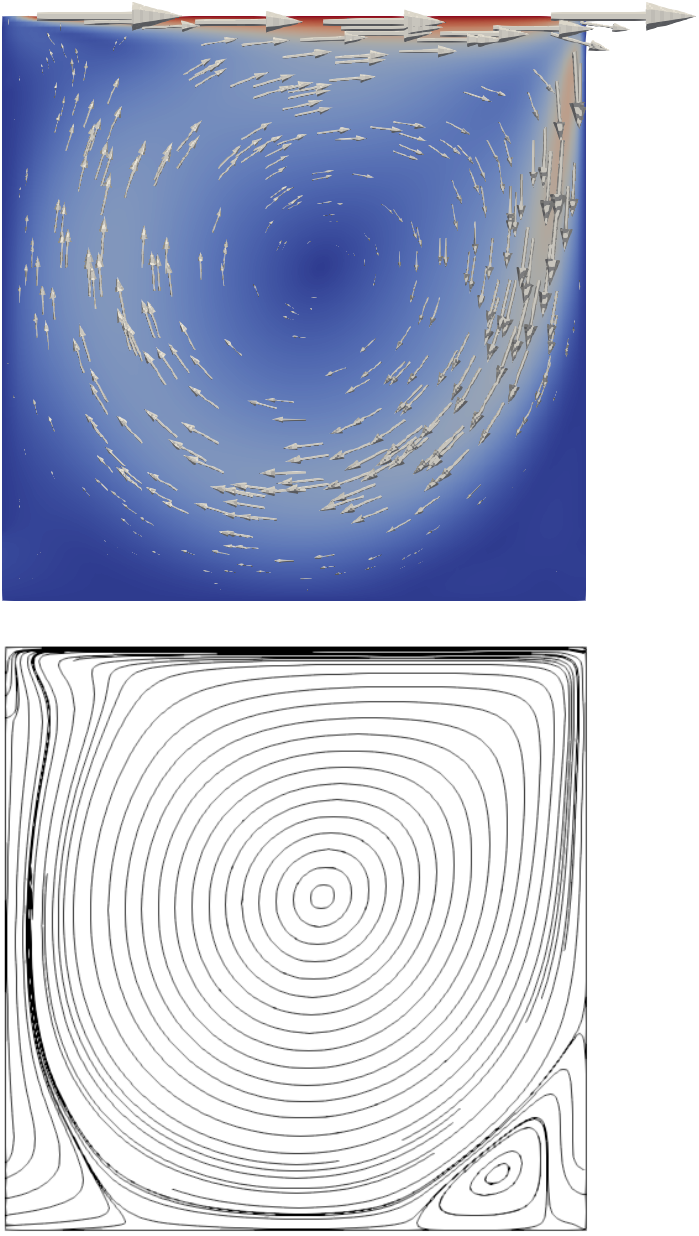}
		\caption{$\mu = 0.05$.}
	\end{subfigure}
	\quad
	\begin{subfigure}[b]{0.305\textwidth}
		\centering
		\includegraphics[width=\textwidth]{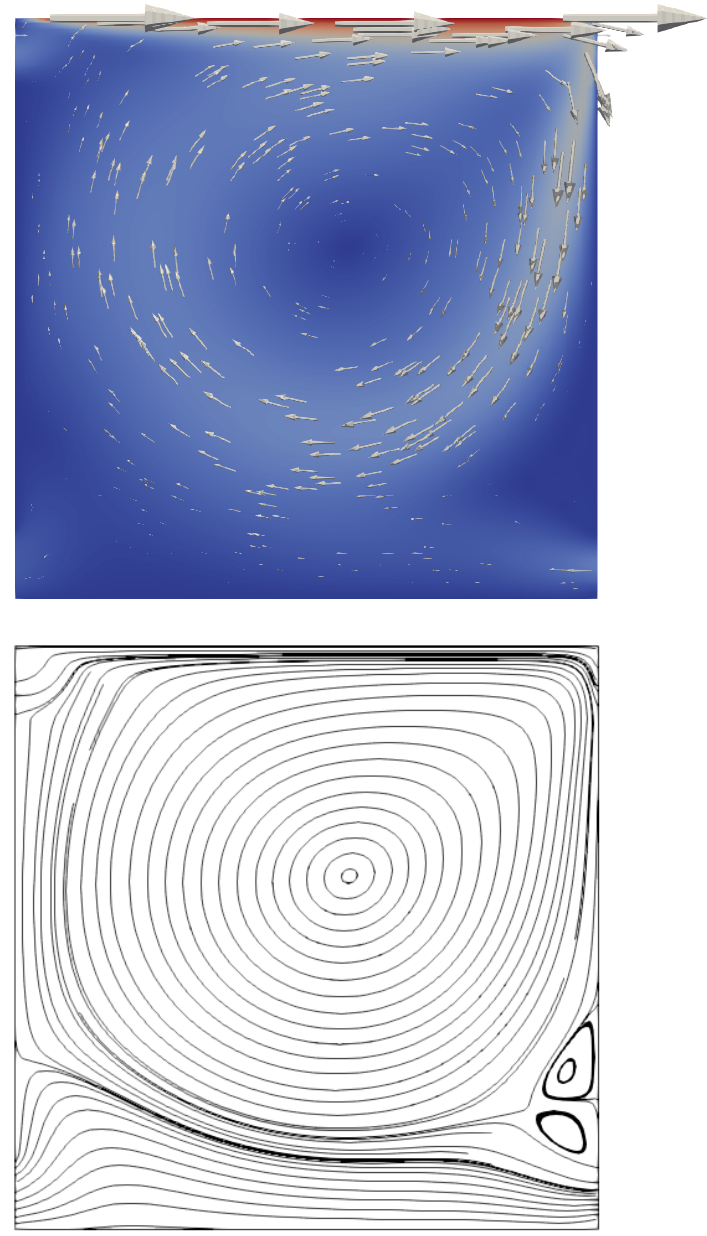}
		\caption{$\mu = 1$.}
	\end{subfigure}
	\caption{\normalsize Velocity field (top) and velocity streamlines (bottom) for the two snapshots in the training set.}
	\label{fig:lid_snaps}
\end{figure}

The training set consists of two snapshots, corresponding to $\mu= 0.05$ and $\mu = 1$. The corresponding velocity fields and streamlines are displayed in Figure \ref{fig:lid_snaps}. 
The POD dimensionality reduction is performed setting the tolerances $\varepsilon_u = 10^{-3}$ and $\varepsilon_p = 0.25$ for velocity and pressure, respectively.
The enrichment procedure with weighting coefficients $\alpha = 0.1, 0.2, \dots, 0.9$ yields an augmented velocity training set containing $11$ elements.

Using all the proposed enriched datasets, the POD constructs reduced spaces of dimension $\nku = 5$ for velocity and $\nkp = 1$ for pressure.
For intermediate values of the parameter $\mu$,  the solution is approximated by the POD-RB using both the standard ROM and the proposed data augmentation strategies for the velocity. 
Table \ref{tab:lidcavity} lists the reduced dimensions and the relative errors, measured in the Euclidean norm, for the velocity and the pressure fields when approximating the solution for $\mu^\star = 0.5$. 
It is straightforward to observe that, in this case, all strategies achieve global errors of order $10^{-2}$ for both velocity and pressure.
However, it is known that the parameters in this problem yield localized variations of the flow, see \cite{TsiolakisPGD2020}.
Hence, a local criterion is introduced,  instead of a global measure of the error, to evaluate the performance of the data augmentation procedures.
More precisely,  the ability of the POD-RB and enriched POD-RB strategies to capture the changes in the flow topology is assessed by comparing the velocity streamlines of the reduced solutions with respect to the high-fidelity ones.
\begin{table}[h!]
	\centering
	\begin{tabular}{ lcccc } 
		\hline
		Data & \multicolumn{2}{c}{Dimension} & \multicolumn{2}{c}{Error} \\[-0.5em]
		augmentation & $\nku$ & $\nkp$ & Velocity & Pressure \\
		\hline 
		None  & $1$ & $1$ & $5.12\cdot10^{-2}$ & $3.15\cdot10^{-2}$    \\ 
		(i)  & $5$ & $1$ & $5.13\cdot10^{-2}$ & $3.54\cdot10^{-2}$    \\ 
		(ii)  & $5$ & $1$ & $3.66\cdot10^{-2}$ & $4.09\cdot10^{-2}$   \\ 
		(iii)  & $5$ & $1$ & $3.38\cdot10^{-2}$ & $2.65\cdot10^{-2}$    \\ 
		\hline
	\end{tabular}
	\caption{Relative errors of the POD-RB approximation, measured in the Euclidean norm, at $\mu^\star = 0.5$ for different strategies of data augmentation: none, standard POD; (i) solenoidal average; (ii) solenoidal average with Oseen enhancement; (iii) linear combination with Oseen enhancement. 
	}
	\label{tab:lidcavity}
\end{table}

\begin{figure}[!h]
	\begin{subfigure}{.3\textwidth}
		\includegraphics[width=\linewidth]{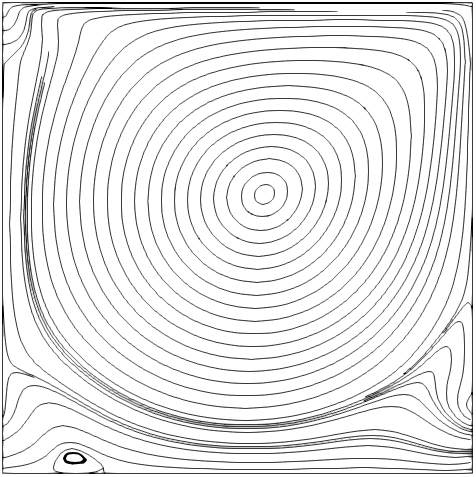}
		\caption{High-fidelity solution.}
		\label{fig:lid_results_05_ref}
	\end{subfigure}\hfill
	\begin{subfigure}{.3\textwidth}
		\includegraphics[width=\linewidth]{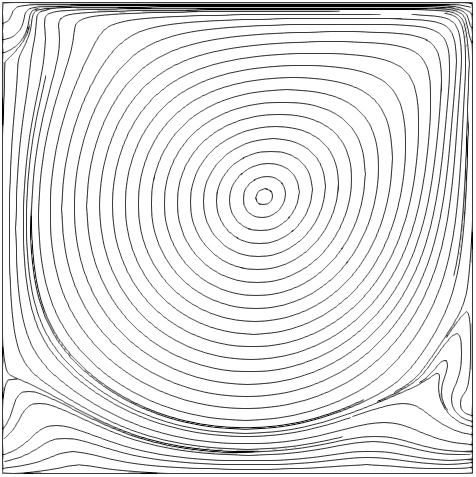}
		\caption{Standard POD-RB.}
		\label{fig:lid_results_05_POD}
	\end{subfigure}\hfill
	\begin{subfigure}{.3\textwidth}
		\includegraphics[width=\linewidth]{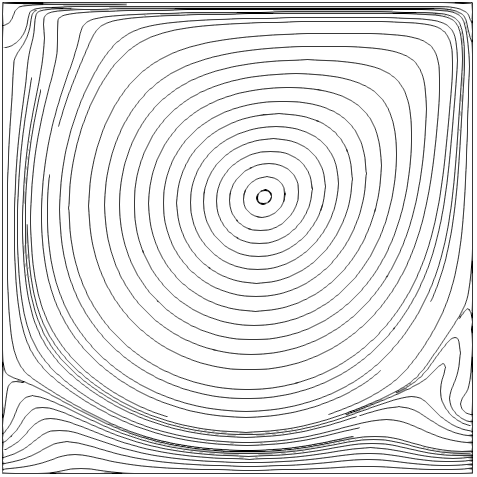}
		\caption{Data augmentation (i).}
		\label{fig:lid_results_05_Sol}
	\end{subfigure}
	\vspace{1mm}

	\begin{minipage}{.3\textwidth}
		$\ $
	\end{minipage} \hfill
	\begin{subfigure}{.3\textwidth}
		\includegraphics[width=\linewidth]{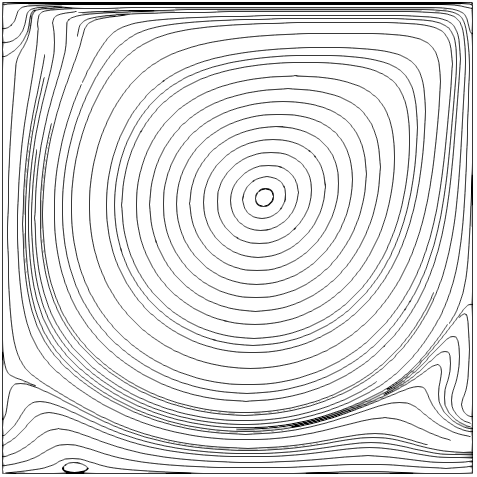}
		\caption{Data augmentation (ii).}
		\label{fig:lid_results_05_Os_ii}
	\end{subfigure}\hfill
	\begin{subfigure}{.3\textwidth}
		\includegraphics[width=\linewidth]{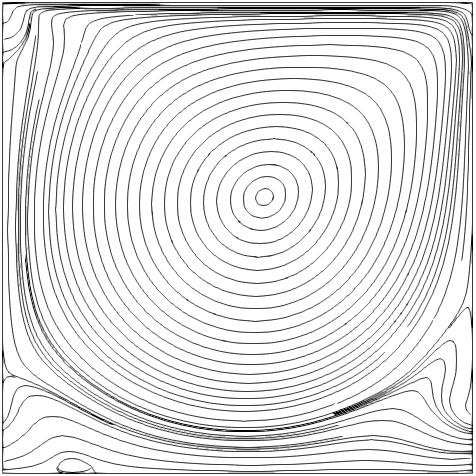}
		\caption{Data augmentation (iii).}
		\label{fig:lid_results_05_Os_iii}
	\end{subfigure}\hfill
	\caption{Velocity streamlines for $\mu^\star = 0.5$ computed using (a) reference high-fidelity solver; (b) POD-RB; POD-RB with data augmentation via (c) solenoidal average; (d) solenoidal average with Oseen enhancement; (e) linear combination with Oseen enhancement. 	
	}
	\label{fig:lid_results_05}
\end{figure}

Figure \ref{fig:lid_results_05} displays the velocity streamlines for $\mu^\star = 0.5$, highlighting the improved ability of the enriched dataset to capture local flow features, such as the nucleation or disappearance of vortices.
In this case, the reference solution in Figure \ref{fig:lid_results_05_ref} shows a vortex in the left-bottom part of the domain. With the considered training set, the standard POD-RB approximation is not able to capture the vortex (Figure \ref{fig:lid_results_05_POD}). Note that including the solenoidal averages to the training set (i.e., enforcing mass conservation in the artificial snapshots) does not provide sufficient additional information (Figure \ref{fig:lid_results_05_Sol}). On the contrary,  Figures \ref {fig:lid_results_05_Os_ii} and \ref{fig:lid_results_05_Os_iii} show that constructing artificial snapshots imposing both mass and momentum conservation principles via the Oseen correction, the reduced-order approximation is capable of identifying the appearence of a vortex where expected.

\begin{figure}[!h]
	\begin{subfigure}{.3\textwidth}
		\includegraphics[width=\linewidth]{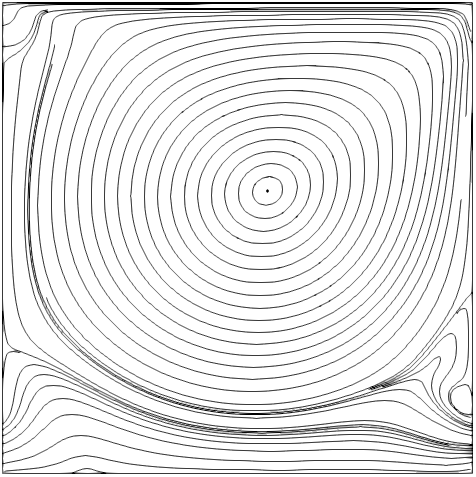}
		\caption{High-fidelity solution.}
		\label{fig:lid_results_07_ref}
	\end{subfigure}\hfill
	\begin{subfigure}{.3\textwidth}
		\includegraphics[width=\linewidth]{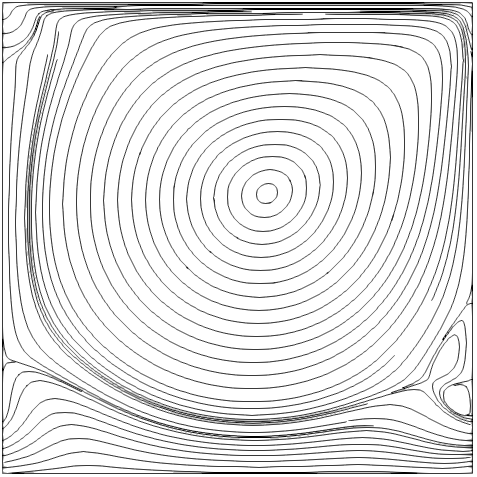}
		\caption{Standard POD-RB.}
		\label{fig:lid_results_07_POD}
	\end{subfigure}\hfill
	\begin{subfigure}{.3\textwidth}
		\includegraphics[width=\linewidth]{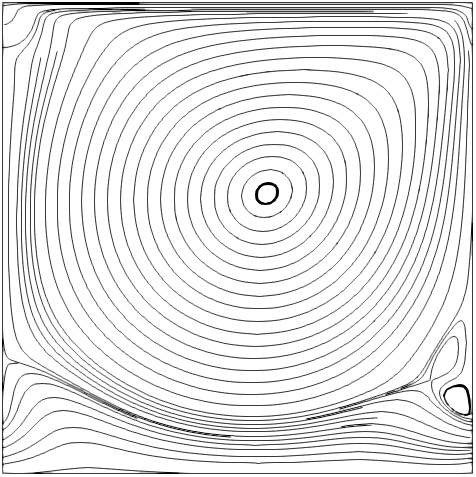}
		\caption{Data augmentation (i).}
		\label{fig:lid_results_07_Sol}
	\end{subfigure}
	\vspace{1mm}

	\begin{minipage}{.3\textwidth}
		$\ $
	\end{minipage} \hfill
	\begin{subfigure}{.3\textwidth}
		\includegraphics[width=\linewidth]{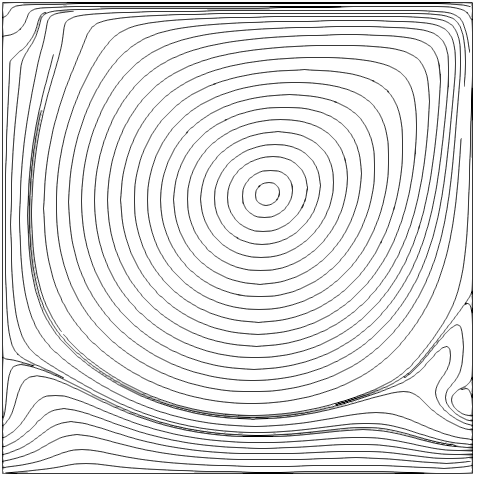}
		\caption{Data augmentation (ii).}
		\label{fig:lid_results_07_Os_ii}
	\end{subfigure}\hfill
	\begin{subfigure}{.3\textwidth}
		\includegraphics[width=\linewidth]{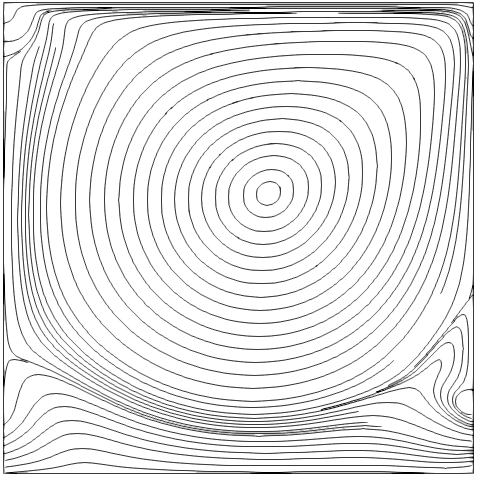}
		\caption{Data augmentation (iii).}
		\label{fig:lid_results_07_Os_iii}
	\end{subfigure}\hfill
	\caption{Velocity streamlines for $\mu^\star = 0.7$ computed using (a) reference high-fidelity solver; (b) POD-RB; POD-RB with data augmentation via (c) solenoidal average; (d) solenoidal average with Oseen enhancement; (e) linear combination with Oseen enhancement. }
	\label{fig:lid_results_07}
\end{figure}

Similarly, the velocity streamlines for $\mu^\star = 0.7$ are shown in Figure \ref{fig:lid_results_07}. In this case,  the POD-RB approximated solutions computed using standard ROM and augmentation with solenoidal averages present two vortices on the bottom right region of the domain (Figures \ref{fig:lid_results_07_POD} and \ref{fig:lid_results_07_Sol}), while only one is visible in the reference solution in Figure \ref{fig:lid_results_07_ref}. The data augmentation relying on both mass and momentum conservation principles introduces in the training set artificial snapshots that improve the overall description of the flow and the reduced approximation displays patterns qualitatively more similar to the high-fidelity solution (Figures \ref {fig:lid_results_07_Os_ii} and \ref{fig:lid_results_07_Os_iii}).
To summarize, incorporating the information about momentum conservation via the Oseen augmentation procedure provides richer artificial snapshots, allowing the resulting enriched POD-RB to capture detailed features of the solution, outperforming the standard POD-RB approach.


\subsection{Parametrized flow past a three-dimensional microswimmer}
\label{sc:microswimmer}

The last example considers a three-dimensional flow past a chiral structure microswimmer, inspired by the studies in \cite{Keaveny-KWS-2013, Giacomini-Sevilla-2020}. 
The microswimmer features a non-axisymmetric helical structure, whence only the data augmentation strategy relying on the Oseen equation with artificial velocity based on the linear combination of snapshots is considered for the comparison with the standard POD-RB.

\begin{figure}[!h]
	\centering
	
	\begin{subfigure}{.45\textwidth}
		\centering
		\includegraphics[height=7cm]{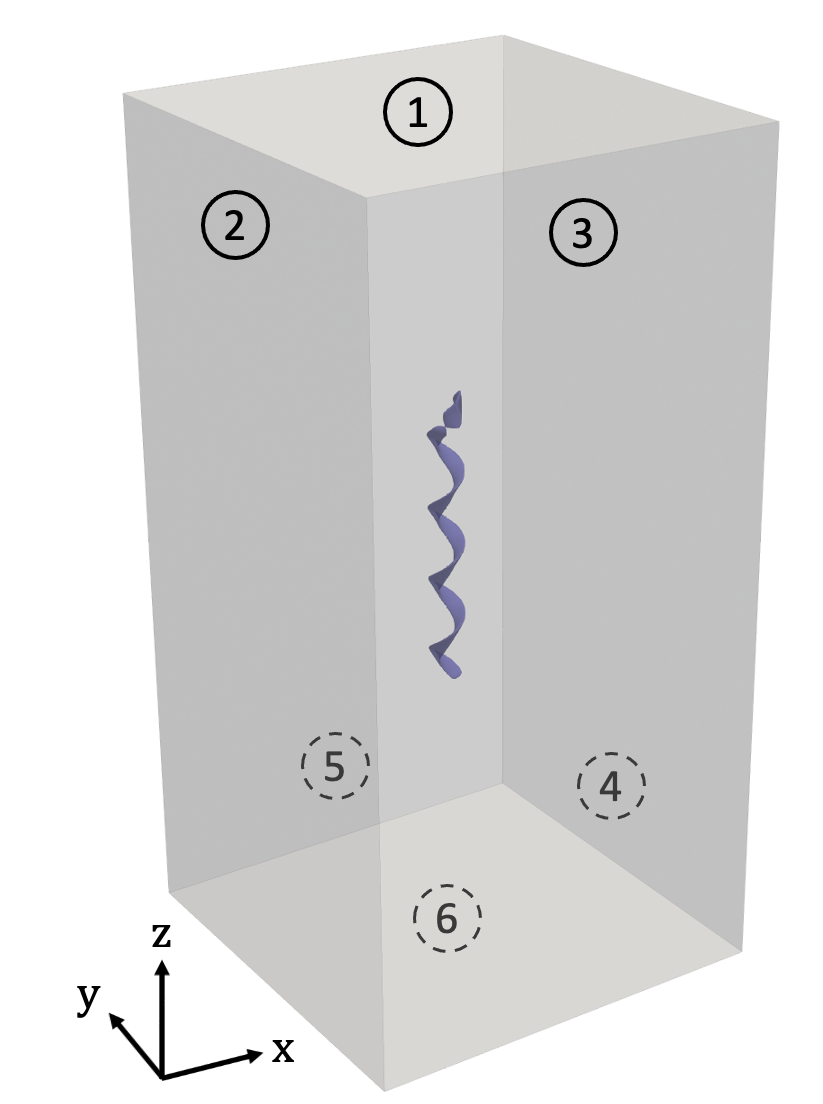}
		\caption{Domain and boundary partitions.}
		\label{fig:swimmer_settingDom}
	\end{subfigure}	
	\begin{subfigure}{.45\textwidth}
		\centering
		\includegraphics[height=7cm]{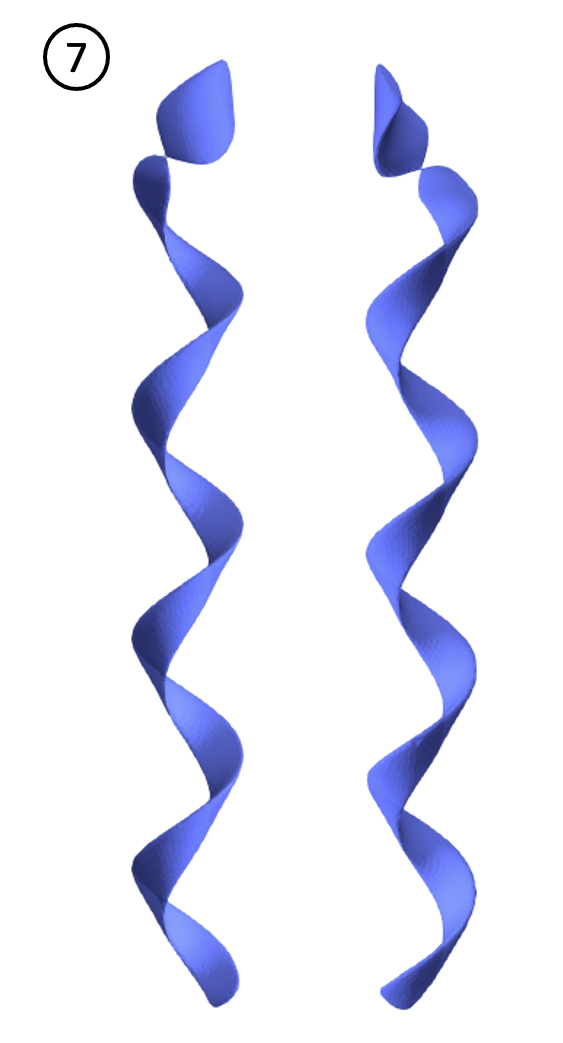}
		\caption{Two views of the geometry of the microswimmer.}
		\label{fig:swimmer_settingGeo}
	\end{subfigure}
	
	\caption{Geometry of the domain and boundary partitions: inlet (faces $1$ and $4$), outlet (faces $2$ and $6$) and free-slip boundaries (faces $3$ and $5$). No slip conditions are imposed on the microswimmer (surface $7$).}
	\label{fig:swimmer_setting}
\end{figure}

The microswimmer, of surface $\mathcal{S}$, length $L = 1.4$ and inner radius $r = 0.057$, is located at the centre of the channel $\mathcal{C} = [-1,1]\times[-1,1]\times[-2,2]$.
The computational domain is therefore the volume $\Omega = \mathcal{C}\setminus\mathcal{S}$. Figure \ref{fig:swimmer_settingGeo} shows the geometry of the microswimmer and the domain, whereas the details of the geometric definition of the microswimmer are presented in \cite{Giacomini-Sevilla-2020}.
\begin{figure}[!h]
	\centering
	\includegraphics[width=.7\textwidth]{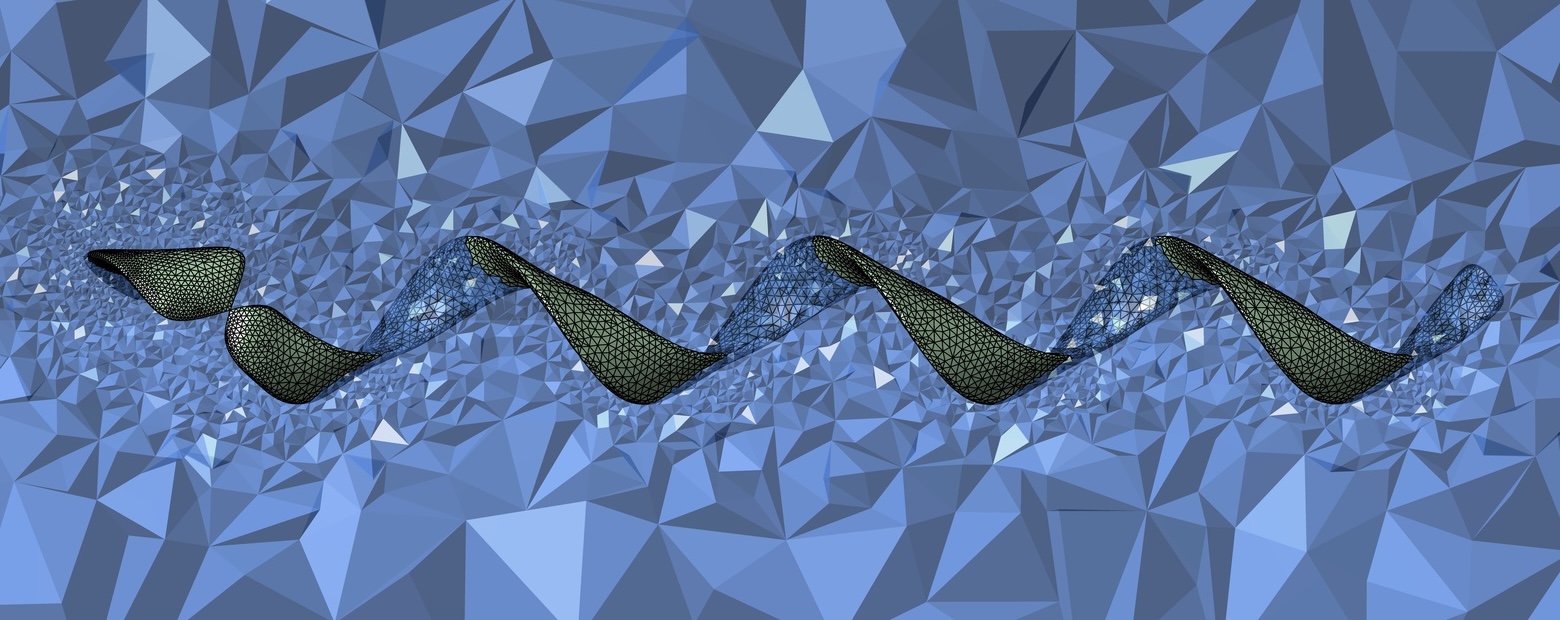}
	\caption{\normalsize Detail of the computational mesh refined in the vicinity of the microswimmer.}
	\label{fig:swimmer_mesh}
\end{figure}	

The flow enters the channel $\mathcal{C}$ through the faces associated with $\{ x = 1 \}$ and $\{ z = 2 \}$. Homogeneous Neumann conditions are imposed on the surfaces of $\mathcal{C}$ at $\{x = -1\}$ and $\{ z = -2\}$ to model the outlet.  Finally, free-slip conditions are enforced on the lateral surfaces of $\mathcal{C}$ at $\{ y = 1 \}$ and $\{ y = -1 \}$, whereas no-slip conditions are imposed on the surface $\mathcal{S}$ of the microswimmer.  A schematic representation of the domain and the boundary conditions is presented in Figure \ref{fig:swimmer_settingDom}.

The flow conditions are parametrized by means of two parameters: the orientation of the inflow velocity and the Reynolds number. 
The former is achieved by setting, at $\{ x = 1 \}$ and $\{ z = 2 \}$, the Dirichlet boundary condition $\bu=\UvecIn(\theta)$ such that
\begin{equation*}
\UvecIn(\theta) = [-\sin\theta, 0, -\cos\theta],
\end{equation*}
where $\theta\in[0,\pi/2]$ is the angle regulating the direction of the flow.  
The Reynolds number is controlled via the kinematic viscosity $\nu$: considering as characteristic velocity the module of $\UvecIn(\theta)$ and as characteristic length $L$,  viscosity values in the interval $[0.14, 14]$ yields a variation of the Reynolds number in the range $[0.1,100]$.

The full-order computations use a $\mathbb{P}_2/\mathbb{P}_1$ finite element pair on a tetrahedral mesh of $463,342$ elements, refined near the microswimmer as displayed in Figure~\ref{fig:swimmer_mesh}. The discreatization leads to $1,903,554$ and $80,976$ degrees of freedom for velocity and pressure, respectively. 
\begin{figure}[h!]
	\centering
	\begin{subfigure}{.32\textwidth}
		\includegraphics[width=\linewidth]{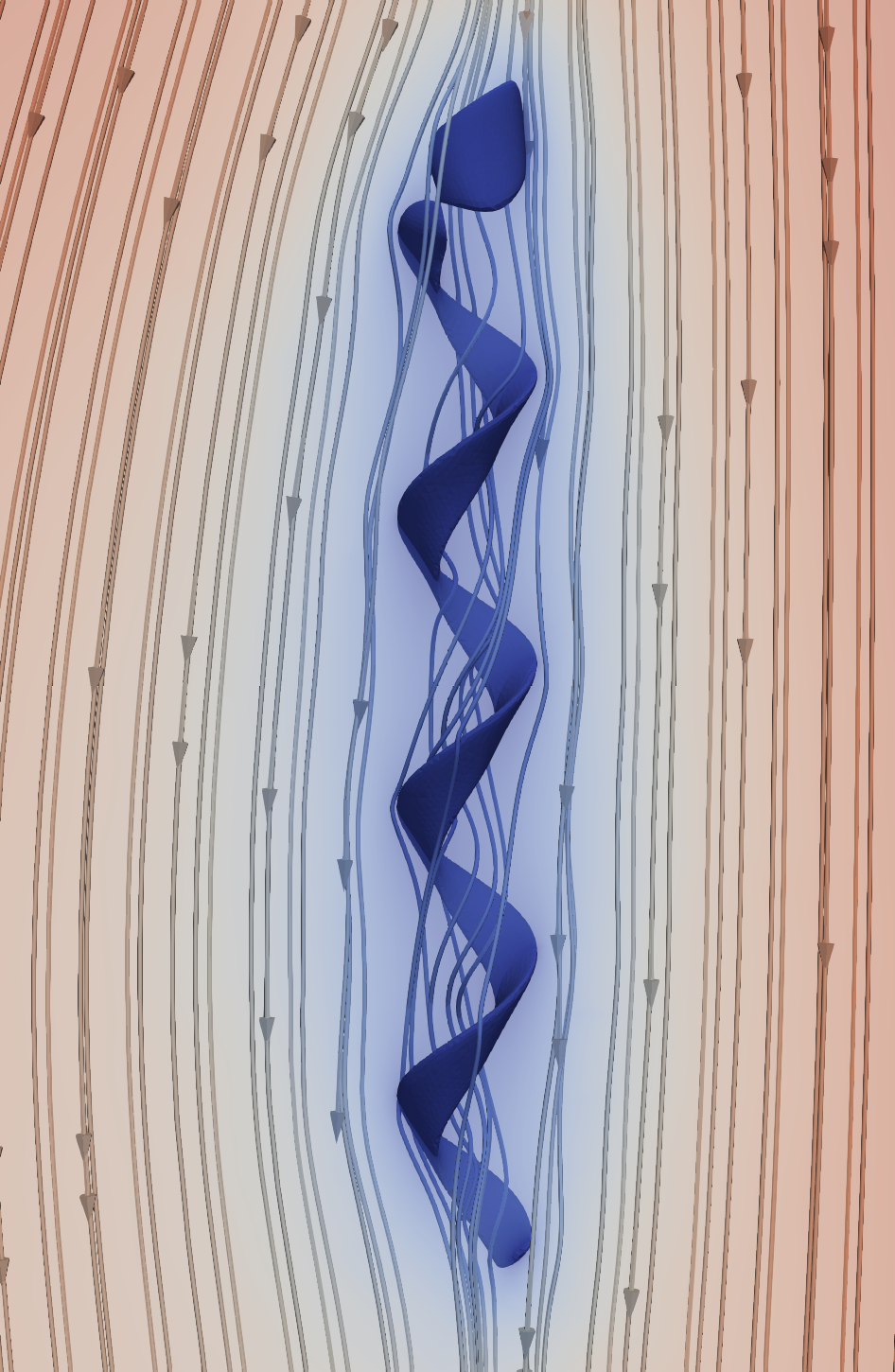}
		\caption{$\theta = 0$, $\Re = 0.1$,.}
	\end{subfigure} \hspace{0.7cm}
	\begin{subfigure}{.32\textwidth}
		\includegraphics[width=\linewidth]{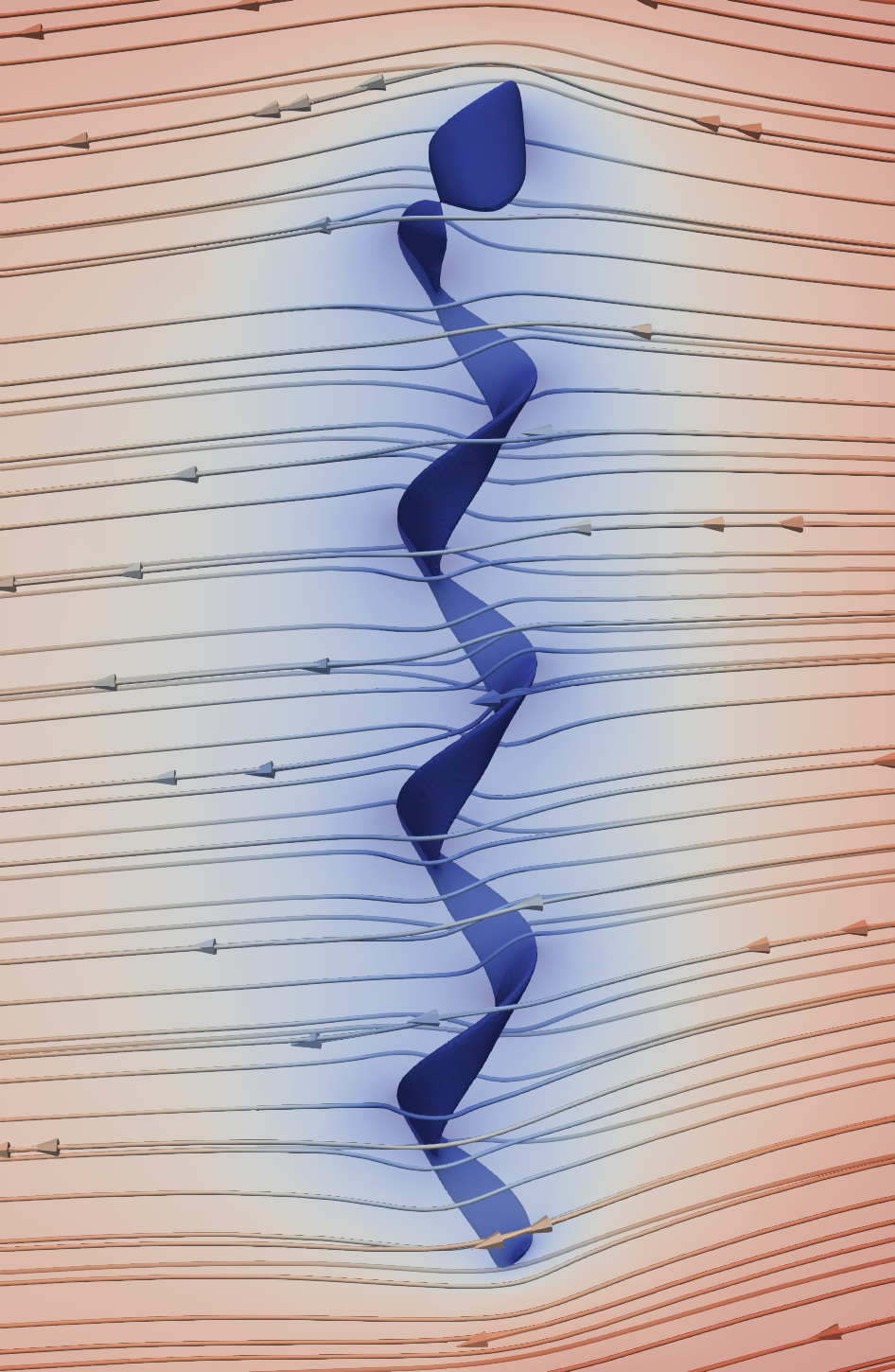}
		\caption{$\theta = \pi/2$, $\Re = 0.1$.}
	\end{subfigure}\hspace{0.5cm}
	\begin{subfigure}{.09\textwidth}
		\hspace{\linewidth}
	\end{subfigure}
	\vspace{1mm}
	
	\begin{subfigure}{.32\textwidth}
		\includegraphics[width=\linewidth]{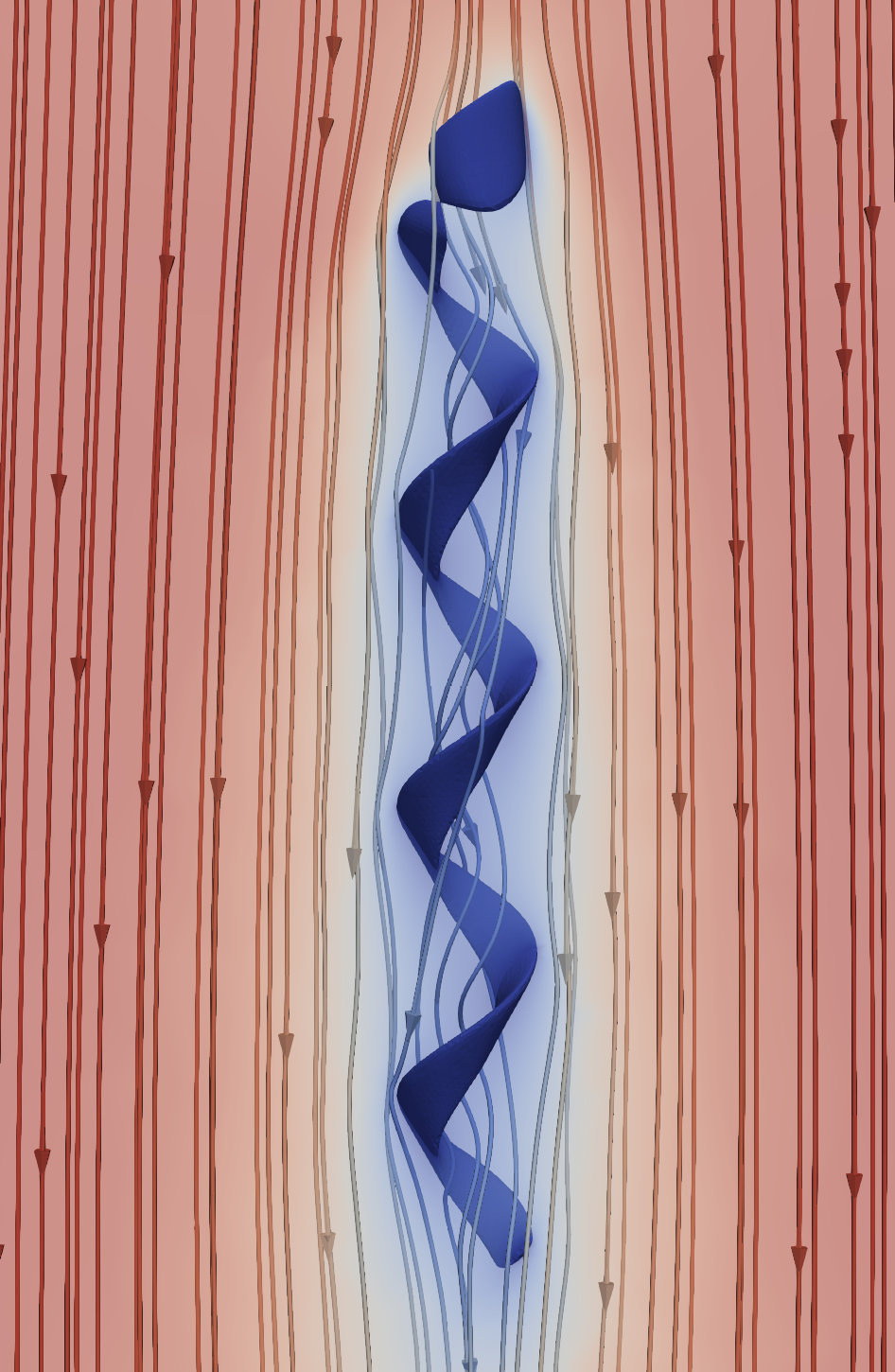}
		\caption{$\theta = 0$, $\Re = 100$.}
	\end{subfigure}\hspace{0.7cm}
	\begin{subfigure}{.32\textwidth}
		\includegraphics[width=\linewidth]{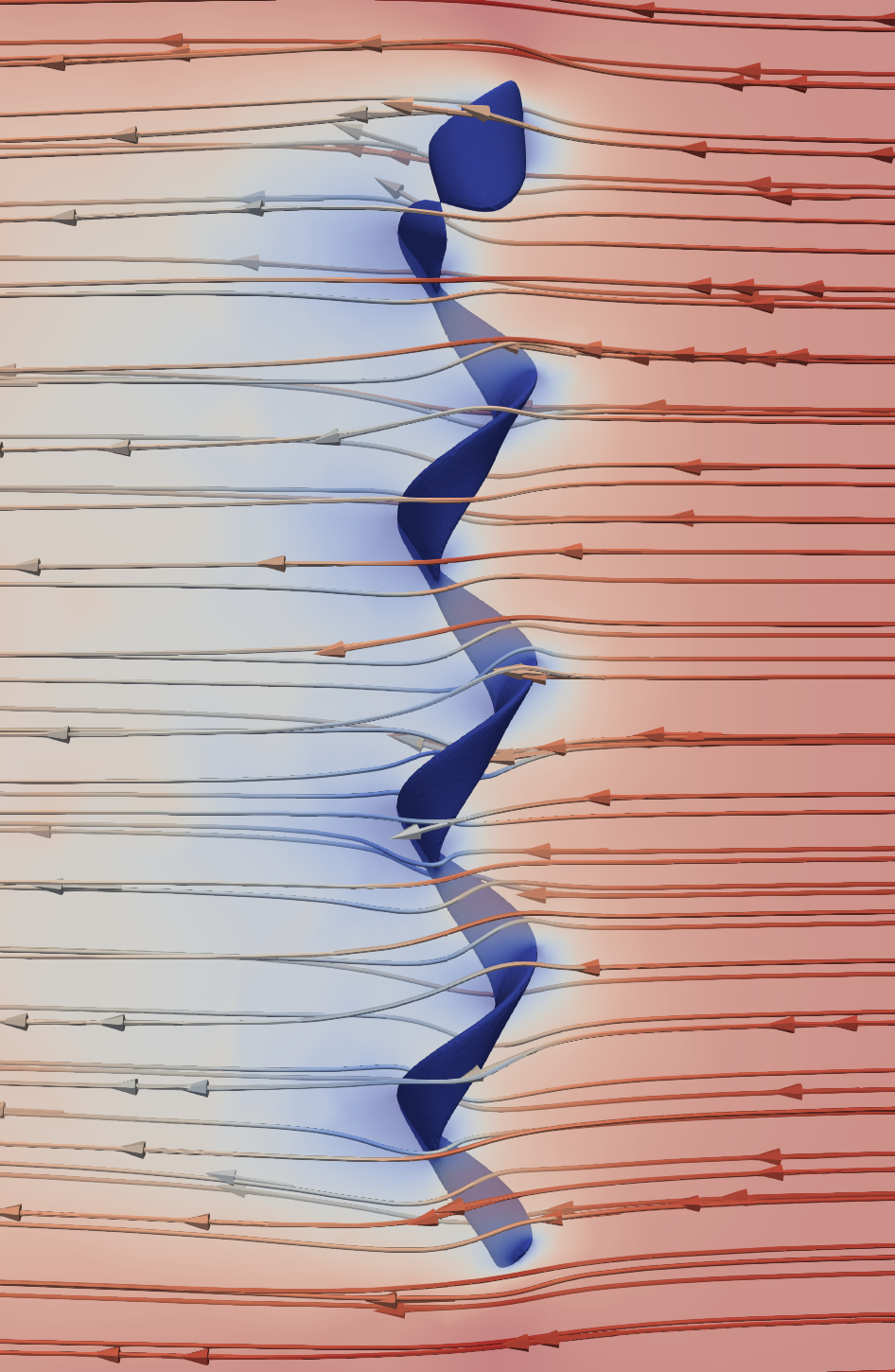}
		\caption{$\theta = \pi/2$, $\Re = 100$.}
	\end{subfigure}\hspace{0.4cm}
	\begin{subfigure}{.09\textwidth}
		\raisebox{0.8cm}{\includegraphics[width=\linewidth]{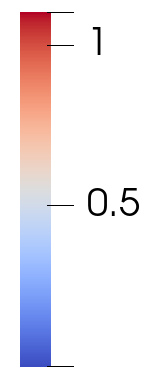}}
	\end{subfigure}
	
	\caption{Streamlines of the flow around the microswimmer for the $4$ snapshots in the training set. Detail at the cross section $\{ y = 0 \}$ for $(x,z)\in[-0.5,0.5]\times[-0.8,0.8]$. The color scale indicates the velocity module.}
	\label{fig:swimmer_family}
\end{figure}

The reduced-order approximations are computed using a training set with four snapshots, corresponding to the combination of the extreme values for each parameter, namely, $(\theta,\Re) \in \{ (0,0.1), (0,100), (\pi/2,0.1),(\pi/2,100) \}$.  The flow field of these four snapshots is depicted in Figure~\ref{fig:swimmer_family}.
Truncation tolerances in the dimensionality reduction are set to $\varepsilon_u = 10^{-3}$ for the velocity training set, and to $\varepsilon_p = 0.4$ for the pressure training set.
As in the previous examples, data augmentation is performed using $9$ linear combinations of the original snapshots, with weighting coefficients $\alpha = 0.1, 0.2, \dots, 0.9$. This leads to an enriched dataset with a total of $40$ snapshots.

The reduced dimensions and the relative errors in the Euclidean norm of the velocity and the postprocessed pressure are listed in Table \ref{tab:swimmer}, when approximating the solution at the intermediate parameters $\theta^\star = \pi/4$ and $\Re^\star = 50$.
The proposed data augmentation strategy enriches the dataset with new, relevant information. In particular, the accuracy of the resulting reduced-order approximation improves, with the error in both velocity and pressure decreasing of almost one order of magnitude with respect to the standard POD-RB one.
\begin{table}[h!]
	\centering
	\begin{tabular}{ lcccc } 
		\hline
		Data & \multicolumn{2}{c}{Dimension} & \multicolumn{2}{c}{Error} \\[-0.5em]
		augmentation & $\nku$ & $\nkp$ & Velocity & Pressure \\
		\hline 
		None (standard POD-RB)  & $3$ & $1$ & $4.20\cdot10^{-1}$ & $5.67\cdot10^{-1}$    \\ 
		Data-enriched POD-RB  & $12$ & $1$ & $5.35\cdot10^{-2}$ & $6.50\cdot10^{-2}$    \\ 
		\hline
	\end{tabular}
	\caption{\normalsize Relative errors of the POD-RB approximation, measured in the Euclidean norm, at $\theta^\star = \pi/4$ and $\Re^\star = 50$.}
	\label{tab:swimmer}
\end{table}

\begin{figure}[h!]
	\begin{subfigure}{.32\textwidth}
		\centering
		\includegraphics[width=\linewidth]{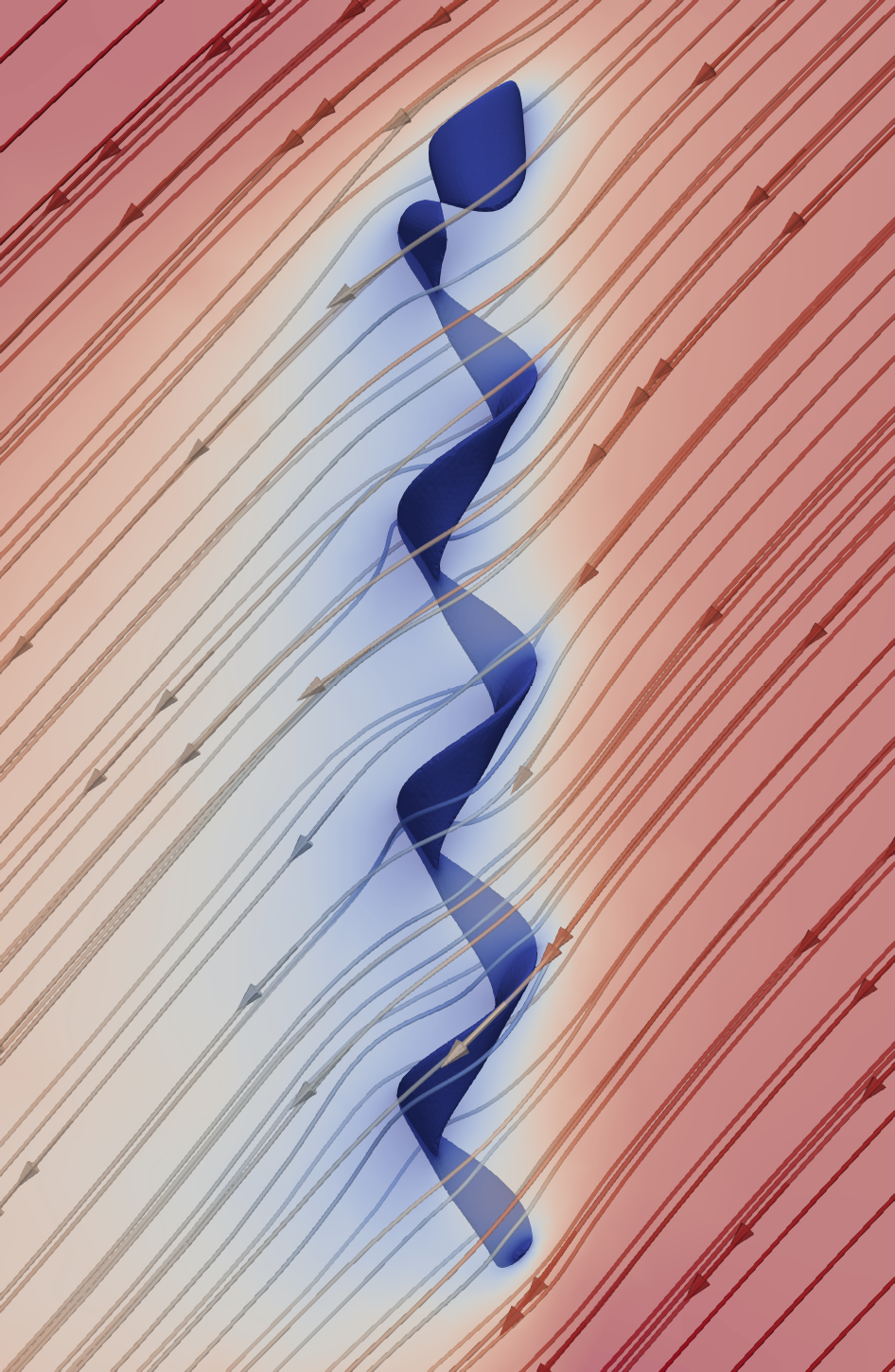}
		\caption{High-fidelity solution.}
		\label{fig:swimmer_sols_i}
	\end{subfigure}\hfill
	\begin{subfigure}{.32\textwidth}
		\centering
		\includegraphics[width=\linewidth]{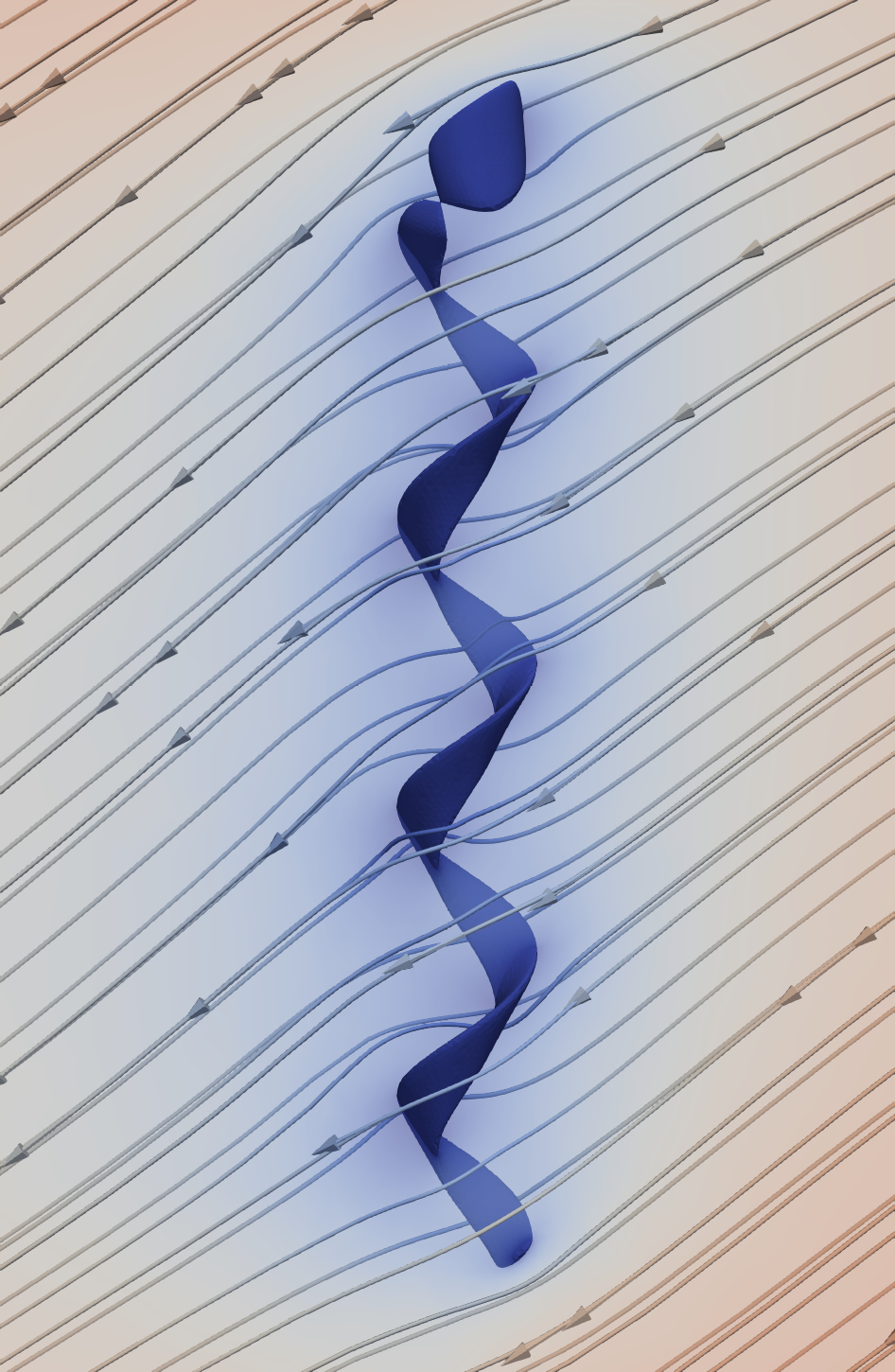}
		\caption{Standard POD.}
		\label{fig:swimmer_sols_ii}
	\end{subfigure}\hfill
	\begin{subfigure}{.32\textwidth}
		\centering
		\includegraphics[width=\linewidth]{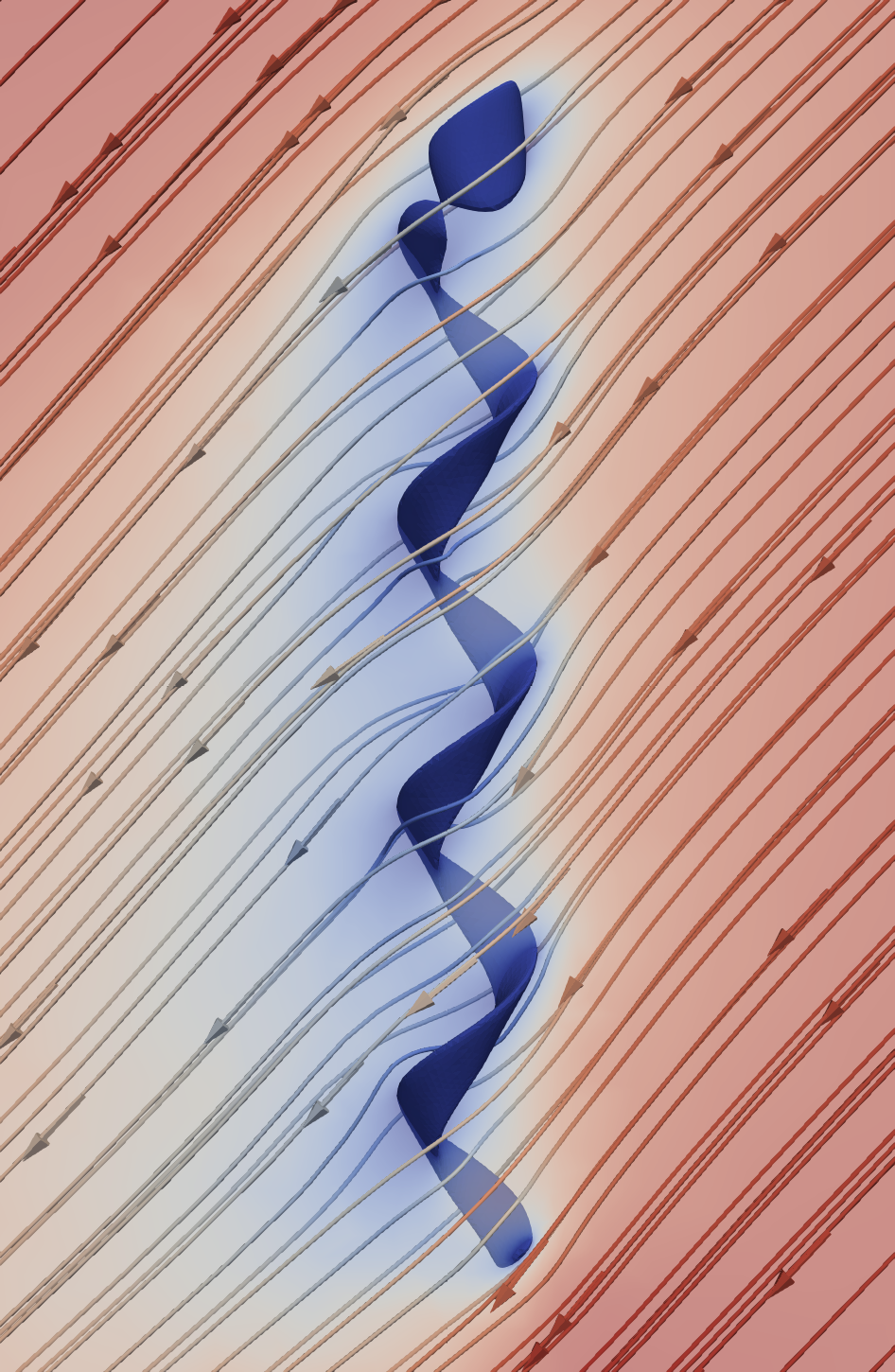}
		\caption{Data-enriched POD.}
		\label{fig:swimmer_sols_iii}
	\end{subfigure}
	\vspace{1mm}

	\caption{Streamlines of the flow around the microswimmer for $\theta^\star = \pi/4$ and $\Re^\star = 50$ using (a) reference high-fidelity solver; the POD-RB approach with (b) standard POD, (c) data-enriched POD. Detail at the cross section $\{ y = 0 \}$ for $(x,z)\in[-0.5,0.5]\times[-0.8,0.8]$.}
	\label{fig:swimmer_sols}
\end{figure}

Finally, the streamlines of the flow around the microswimmer for $\theta^\star = \pi/4$ and $\Re^\star = 50$ are reported in Figure \ref{fig:swimmer_sols}. The results confirm that the data augmentation strategy introduces new information in the reduced space. The corresponding data-enriched POD-RB solution (Fig. \ref{fig:swimmer_sols_iii}) outperforms the standard POD-RB (Fig. \ref{fig:swimmer_sols_ii}) in approximating the full-order finite element solution in Figure \ref{fig:swimmer_sols_i}.

\section{Discussion on accuracy and efficiency of the enriched POD-RB}
\label{sec:Discussion}

The described enriched POD-RB relying on physics-based data augmentation procedures aims to mimick the accuracy of a standard POD-RB approximation constructed starting from a richer dataset.  Indeed, the goal of data augmentation is to construct new snapshots that likely approximate the solution at intermediate parametric values, at a fraction of the cost of solving an additional full-order problem. Of course, high-fidelity data provided by a full-order solver are to be preferred to any artificial snapshot.
Nonetheless, the availability of such additional information is subject to (i) the extra computational cost for the calculation of new high-fidelity snapshots, and (ii) tailored sampling algorithms to select the parameters for which new evaluations of the full-order solver are required \citep{Willcox-BWG-08}.
In the following, uniform samplings are considered for both standard and enriched POD.

In this context, the target accuracy that the enriched POD-RB can expect to achieve is bounded by the accuracy of the standard POD-RB computed with a larger number of snapshots.  
For the 3D case in Section \ref{sc:microswimmer}, the enriched POD-RB with different numbers of linear combinations of the original snapshots is compared with the standard POD-RB with $4$ and $24$ uniformly sampled snapshots. 
For the enriched POD, $9$, $5$ and $1$ linear combinations of the $4$ original high-fidelity snapshots are considered, leading to enriched datasets with a total of $40$, $24$ and $8$ snapshots, respectively.
The results are reported in Table \ref{tab:swimmerComb}: besides the conclusions already presented in Section  \ref{sc:microswimmer} on the ability of the enriched POD-RB with $9$ combinations to lower by almost one order of magnitude the error of the standard POD-RB with the same amount of high-fidelity snapshots,  it is worth noticing that the gain in accuracy is preserved even when fewer artificial snapshots are created. 
Indeed, while decreasing the number of combinations considered increases the error in velocity from $5.35\cdot10^{-2}$ to $7.04\cdot10^{-2}$ and the error in pressure from $6.50\cdot10^{-2}$ to $8.55\cdot10^{-2}$, they all remain below the corresponding errors of the standard POD-RB at $4.20\cdot10^{-1}$ and $5.67\cdot10^{-1}$ , respectively. 
Hence, by solely solving $4$ additional linear systems \eqref{eq:discreteOseen}, the data augmentation provides significantly more accurate solutions. 
\begin{table}[h!]
	\centering
	\begin{tabular}{ lcccc } 
		\hline
		& \multicolumn{2}{c}{Dimension} & \multicolumn{2}{c}{Error} \\[-0.5em]
		& $\nku$ & $\nkp$ & Velocity & Pressure \\
		\hline 
		\multicolumn{5}{c}{Standard POD-RB} \\
		\hline
		24 snapshots  & $15$ & $1$ & $1.75\cdot10^{-2}$ & $1.64\cdot10^{-2}$    \\ 
		4 snapshots  & $3$ & $1$ & $4.20\cdot10^{-1}$ & $5.67\cdot10^{-1}$    \\ 
		\hline 
		\multicolumn{5}{c}{Enriched POD-RB using 4 snapshots} \\
		\hline 
		9 combinations  & $12$ & $1$ & $5.35\cdot10^{-2}$ & $6.50\cdot10^{-2}$    \\ 
		5 combinations  & $10$ & $1$ & $6.88\cdot10^{-2}$ & $9.52\cdot10^{-2}$   \\ 
		1 combination  & $5$ & $1$ & $7.04\cdot10^{-2}$ & $8.55\cdot10^{-2}$    \\ 
		\hline
	\end{tabular}
	\caption{Relative errors of the POD-RB approximation, measured in the Euclidean norm, at $\theta^\star = \pi/4$ and $\Re^\star = 50$ for the standard POD-RB with different number of snapshots and the enriched POD-RB with different number of combinations to generate artificial snapshots.
	}
	\label{tab:swimmerComb}
\end{table}

The standard POD-RB approximation with $24$ snapshots is now compared with the enriched POD-RB computed with $4$ original snapshots and $5$ combinations to generate artificial snapshots.  The resulting training sets have the same number of entries.  
As expected, standard POD-RB with $24$ high-fidelity snapshots leads to a more accurate approximation, achieving errors of $1.75\cdot10^{-2}$ in velocity and $1.64\cdot10^{-2}$ in pressure, whereas the enriched POD-RB yields errors of $6.88\cdot10^{-2}$ and $9.52\cdot10^{-2}$ in velocity and pressure, respectively.
Nonetheless, for the standard POD-RB, the computation of the additional $20$ high-fidelity snapshots solving equation \eqref{eq:NonLinearSystem}, with either $3$ or $4$ iterations of the Newton-Raphson method per snapshot, entails between $60$ and $80$ extra linear systems to be solved.
On the contrary,  the data augmentation procedure requires the solution of $20$ additional linear systems \eqref{eq:discreteOseen}, yielding a reduction of the computational effort between $67\%$ and $75\%$.

\begin{remark}
	For the above assessment of the cost of the data augmentation procedure, only the solution of the Navier-Stokes equations, for the high-fidelity snapshots, and of the Oseen equations, for the artificial ones, are considered.
	Indeed, the remaining steps involved in the construction of the artificial snapshots rely on trivial operations (e.g., vector-scalar product, sum of vectors, ...) whose computational cost is negligible with respect to the solution of large-scale linear systems of equations.
	Similarly, the solenoidal average only entails the solution of the linear scalar elliptic partial differential equation \eqref{eq:SolPoisson}. This step introduces a limited number of additional unknowns and leads to a symmetric positive definite matrix that can be easily factored to accelerate the solution of the data augmentation problems, yielding a negligible overhead cost.
\end{remark}

\section{Concluding remarks}
\label{sec:Conclusion}

The methodologies presented in this paper aim at producing, with limited additional computational cost, new, artificial snapshots to enrich the training set for a posteriori ROMs.
The ideas behind the proposed strategies are inspired in the physical knowledge of the solution,  that is, enforcing mass and momentum conservation principles.

Two approaches are proposed, the one referred to as solenoidal geometric average constructs artificial snapshots fulfilling mass conservation equation,  whereas the one labeled Oseen enhancement enforces both mass and momentum conservation. 
Despite the apparent soundness of the artificial snapshots produced by the solenoidal geometric average, it was found that, in the context of steady-state incompressible Navier-Stokes, they do not bring to the reduced basis gainful features. That is, no more than the linear combinations proposed by the standard POD-RB computed using the original snapshots. 
On the contrary, accounting for both mass and momentum conservation principles via the Oseen enhancement is actually introducing in the training set new, relevant information, particularly useful to discover the features of parametric configurations not seen during the training procedure.

This opens the way to create data augmentation strategies, tailored for each specific application, to boost POD-RB and related techniques. This is expected to release these ROMs from the computational burden of populating a training set with numerous full-order solutions.
Of course, a rigorous analysis to identify the emerging approximation properties obtained by adding new, artificial snapshots to the dataset will provide a solid background to the discussed physics-informed data augmentation strategies.
Moreover,  further studies involving problems with multiple parameters are to be performed in order to assess the scalability of the proposed enrichment strategies to high-dimensional parametric spaces.
	
Finally, the suitability and performance of the proposed methodologies in the context of transient and turbulent flows need to be throroughly studied. On the one hand, data augmentation techniques are expected to yield significant computational gains in time-dependent parametric problems, where the augmented data is reused many times. On the other hand, the individual relevance of mass and momentum conservation principles are to be assessed when dynamic effects and multi-scale phenomena are involved.

\section*{Acknowledgments}
This work was partially supported by the Marie Sklodowska-Curie Actions (Doctoral Network with Grant agreement No. 101120556 to AM, SZ and PD), the Spanish Ministry of Science, Innovation and Universities and the Spanish State Research Agency MICIU/AEI/10.13039/501100011033 (Grant agreement No. PID2020-113463RB-C32 to SZ; PID2020-113463RB-C33 to AM, MG and PD; CEX2018-000797-S to SZ, MG and PD; TED2021-132021B-I00 to MG) and the Generalitat de Catalunya (Grant agreement No. 2021-SGR-01049).
AM acknowledges the support of the Spanish Ministry of Science, Innovation and Universities through the Margarita Salas fellowship.
AM and MG are Fellows of the Serra H\'unter Programme of the Generalitat de Catalunya.

\section*{Conflict of interest}

The authors declare no potential conflict of interests.

\bibliographystyle{elsarticle-harv}
\bibliography{Ref}

\end{document}